\documentclass[11pt]{article}
\usepackage{comment}
\usepackage{enumitem}
\usepackage{tikz}
\usetikzlibrary{shapes}
\usetikzlibrary{snakes}
\usepackage{hyperref}
\usepackage{microtype,a4wide}
\usepackage{boxedminipage}
\usepackage{xspace}
\hypersetup{
    pdftitle=   {Computing Small Pivot-Minors},
   pdfauthor=  {Konrad K. Dabrowski, Fran\c{c}ois Dross, Jisu Jeong, Mamadou Moustapha Kant\'e, 
O-joung~Kwon, Sang-il Oum, Dani\"el Paulusma}
}
\usepackage{todonotes}
\usepackage{authblk}
\usepackage{amsthm,amssymb,amsmath}
\usepackage{blkarray}
\usepackage{mathabx}
\newtheorem{theorem}{Theorem}[section]
\newtheorem{lemma}[theorem]{Lemma}
\newtheorem{proposition}[theorem]{Proposition}

\newcommand\abs[1]{\lvert #1\rvert}

\newcommand{\problemdef}[3]{
	\begin{center}
		\begin{boxedminipage}{.82\textwidth}
			\textsc{{#1}}\\[1pt]
			\begin{tabular}{ r p{0.7\textwidth}}
				\textit{~~~~Instance:} & {#2}\\
				\textit{Question:} & {#3}
			\end{tabular}
		\end{boxedminipage}
	\end{center}
}

\def\bw_#1{{\overline{BW_{#1}}}}

\newcommand{\ssi}{{\subseteq_i}}
\newcommand{\NP}{{\sf NP}}
\newcommand\Bin{\operatorname{Bin}}

\newcommand\pivot{\wedge}
\def\bw_#1{{\overline{BW_{#1}}}}

\begin{document}
\title{Computing Pivot-Minors\thanks{An extended abstract of this paper appeared in the proceedings of WG 2018~\cite{Dabrowski2018}.
Dabrowski and Paulusma 
were supported by the Leverhulme Trust (RPG-2016-258). This work was mainly done when Jeong was in KAIST.  Kant\'e was supported by the French Agency for Research under the projects DEMOGRAPH (ANR-16-CE40-0028) and ASSK (ANR-18-CE40-0025). Kwon was supported by the National Research Foundation of Korea funded by the Ministry of Science and ICT (No. NRF-2021K2A9A2A11101617 and No. RS-2023-00211670). Kwon and Oum were supported by the Institute for Basic Science (IBS-R029-C1).
Dross was supported by the European Research Council (ERC) under the European Union's Horizon 2020 research and innovation programme Grant Agreement 714704.}}

\author[1]{Konrad K. Dabrowski}
\affil[1]{School of Computing, Newcastle University, UK}

\author[2]{Fran\c{c}ois Dross}
\affil[2]{LaBRI, Universit\'e de Bordeaux, France}

\author[3]{Jisu Jeong}
\affil[3]{Clova AI Research, NAVER Corp, Seongnam,~Korea}

\author[4]{Mamadou Moustapha Kant\'e} 
\affil[4]{Universit\'e Clermont Auvergne, LIMOS, CNRS, Aubi\`ere, France}

\author[5,6]{O-joung Kwon}
\affil[5]{Department of Mathematics, Hanyang University, Seoul,~Korea}

\author[6,7]{Sang-il Oum}
\affil[6]{Discrete Mathematics Group,
  Institute for Basic Science (IBS), Daejeon,~Korea}
\affil[7]{Department of Mathematical Sciences, KAIST, Daejeon,~Korea}
\author[8]{Dani\"el Paulusma}
\affil[8]{Department of Computer Science, Durham University, UK}

\renewcommand\footnotemark{}
\date\today
\maketitle

  \setcounter{footnote}{1}%
 \footnotetext{E-mail addresses: 
    \texttt{konrad.dabrowski@newcastle.ac.uk} (Dabrowski),
    \texttt{francois.dross@u-bordeaux.fr} (Dross),
    \texttt{jisujeong89@gmail.com} (Jeong),
    \texttt{mamadou.kante@uca.fr} (Kant\'e),
    \texttt{ojoungkwon@hanyang.ac.kr} (Kwon),
    \texttt{sangil@ibs.re.kr} (Oum),
    \texttt{daniel.paulusma@durham.ac.uk} (Paulusma).
    }

\begin{abstract}
\noindent
A graph~$G$ contains a graph~$H$ as a pivot-minor if~$H$ can be obtained from~$G$ by applying a sequence of vertex deletions and edge pivots. Pivot-minors play an important role in the study of rank-width.
Pivot-minors have mainly been studied from a structural perspective. 
In this paper we perform the first systematic computational complexity study of pivot-minors.
We first prove that the {\sc Pivot-Minor} problem, which asks if a given graph~$G$ contains a pivot-minor isomorphic to a given graph~$H$, is \NP-complete. If~$H$ is not part of the input, we denote the problem by $H$-{\sc Pivot-Minor}.
We give a certifying polynomial-time algorithm for {\sc $H$-Pivot-Minor} when
\begin{itemize}
    \item $H$ is an induced subgraph of $P_3+tP_1$ for some integer $t\geq 0$,
    \item $H=K_{1,t}$ for some integer $t\geq 1$, or
    \item $|V(H)|\leq 4$ except when $H \in \{K_4,C_3+\nobreak P_1\}$.
\end{itemize}
Let~${\cal F}_H$ be the set of induced-subgraph-minimal graphs that contain a pivot-minor isomorphic to $H$.
To prove the above statement, we either show that there is an integer $c_H$ such that all graphs in ${\cal F}_H$ have at most $c_H$ vertices, or we  determine ${\cal F}_H$ precisely, for each of the above cases.

\end{abstract}

\section{Introduction}\label{s-intro}

Computing whether a graph~$H$ appears as a ``pattern'' inside some other graph~$G$ is a well-studied problem in the area of structural and algorithmic graph theory.
The definition of a pattern depends on the set of graph operations that we are allowed to use.
For instance, if we can obtain~$H$ from~$G$ via a sequence of vertex deletions, edge deletions and edge contractions, then~$G$ contains~$H$ as a minor.
The {\sc Minor} problem is that of testing whether a given graph~$G$ contains a minor isomorphic to a given graph~$H$.
This problem is known to be \NP-complete even if~$G$ and~$H$ are trees of small diameter~\cite{MT92}.
Hence, it is natural to fix the graph~$H$, and let the input consist of only~$G$.
This leads to the {\sc $H$-Minor} problem. A celebrated result of Robertson and Seymour~\cite{RS95} states that $H$-{\sc Minor} can be solved in cubic time for every graph~$H$.
If we only allow vertex deletions and edge contractions, then we obtain the $H$-{\sc Induced Minor} problem.
In contrast to {\sc $H$-Minor}, the {\sc $H$-Induced Minor} problem is \NP-complete for some graphs~$H$; see~\cite{FKMP95} for the smallest known``hard'' graph~$H$, which has $68$ vertices.
Other well-known containment relations include containing a graph~$H$ as a contraction, an induced subgraph, a subdivision, or
an (induced) topological minor; see, for example,~\cite{BV87,GKMW11,LLMT09,LMT12,HKPST12} for a number of complexity results for these relations.

\subsection{Our Focus}

We consider the pivot-minor containment relation. This relation is defined as follows.
The \emph{local complementation} at a vertex~$u$ in a graph~$G$ replaces every edge of the subgraph induced by the set $N_G(u)$ of neighbours of~$u$ by a non-edge, and vice versa.
We denote the resulting graph by~$G*u$.
An \emph{edge pivot} is the operation that takes an edge~$uv$, first applies a local complementation at~$u$, then at~$v$, and then at~$u$ again.
We denote the resulting graph by 
$$G\wedge uv=((G*u)*v)*u.$$ It is known and can be readily checked that $((G*u)*v)*u=((G*v)*u)*v$~\cite{Ou05}.

\begin{figure}
  \begin{center}
  \begin{tikzpicture}
    \tikzstyle{v}=[circle,draw,fill=black!50,inner sep=0pt,minimum width=4pt]
    \draw (0,0)node[v](v){};
    \foreach \i in {1,2,3,4,6,7} {
      \draw (90-360/14+360*\i/7:1) node [v] (v\i){};
      }
    \draw (0,-0.5) node [v] (v5){};
      \node at (v7) [label=above:$v$]{};
      \node at (v1) [label=above:$u$]{};
      \draw (v3)--(v2)--(v1)--(v7)--(v6)--(v4);
      \foreach \i in {7,1,2} {
      \draw (v5)--(v\i);
    }
    \draw (v7)--(v)--(v1);
    \draw (v1)--(v3)--(v6);
    \node at (v4) [label=below:$G$]{};
  \end{tikzpicture}
  $\quad\quad$
  \begin{tikzpicture}
    \tikzstyle{v}=[circle,draw,fill=black!50,inner sep=0pt,minimum width=4pt]
    \draw (0,0)node[v](v){};
    \foreach \i in {1,2,3,4,6,7} {
      \draw (90-360/14+360*\i/7:1) node [v] (v\i){};
      }
    \draw (0,-0.5) node [v] (v5){};
      \node at (v7) [label=above:$u$]{};
      \node at (v1) [label=above:$v$]{};
      \draw (v3)--(v2)--(v1)--(v7)--(v6)--(v4);
      \foreach \i in {7,6,1,3} {
      \draw (v5)--(v\i);
    }
    \draw (v1)--(v3);
    \draw (v6)--(v2);
    \draw (v7)--(v)--(v1);
    \foreach \i in {2,3,6}{
      \draw (v\i)--(v);
    }
    \node at (v4) [label=below:$G\pivot uv$]{};
  \end{tikzpicture}
  \end{center}
\caption{An example of a graph before (left) and after (right) pivoting the edge $uv$.}
\label{f-pivotexample}
\end{figure}
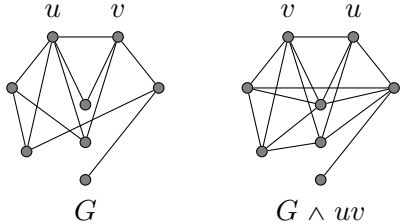

Alternatively, we can define the edge pivot operation as follows.
Let $S_u$ be the set of neighbours of~$u$ in $V(G)\setminus \{v\}$ that are not adjacent to~$v$.
Let $S_v$ be the set of neighbours of~$v$ in $V(G)\setminus \{u\}$ that are not adjacent to~$u$.
Finally, let $S_{uv}$  be the set of common neighbours of~$u$ and~$v$.
We replace every edge between any two vertices in distinct sets from $\{S_u,S_v,S_{uv}\}$ by a non-edge and vice versa.
Next, we delete every edge between~$u$ and~$S_u$ and add every edge between~$u$ and~$S_v$.
Similarly, we delete every edge between~$v$ and~$S_v$ and add every edge between~$v$ and~$S_u$.
See \figurename~\ref{f-pivotexample} for an example of an edge pivot.

A graph~$G$ contains a graph~$H$ as a \emph{pivot-minor} if~$G$ can be modified into $H$ by a sequence of vertex deletions and edge pivots. If $G$ does not contain a pivot-minor isomorphic to $H$, then $G$ is \emph{$H$-pivot-minor-free}.

\subsection{Related Work}

Pivot-minors were initially
called $p$-reductions by Bouchet~\cite{Bo94}. They have 
been studied from a structural perspective, as they form a very suitable tool for working with rank-width~\cite{Ou05,OS06}.
Rank-width is a well-known width parameter (see~\cite{Ou16} for a survey). Pivot-minors play a similar role for rank-width as minors do for treewidth.
Oum~\cite{Ou08} showed that for every positive constant~$k$, the class of graphs of rank-width at most~$k$ is well-quasi-ordered under the pivot-minor relation.
Kwon and Oum~\cite{KO14} proved that every graph of rank-width at most~$k$ is a pivot-minor of a graph of treewidth at most~$2k$, and that a graph of linear rank-width at most~$k$ is a pivot-minor of a graph of path-width at most~$k+1$.

Pivot-minors are closely related to so-called vertex-minors, introduced in the nineties as $\ell$-reductions by Bouchet~\cite{Bo94}.
A graph~$G$ contains a graph~$H$ as a \emph{vertex-minor} if~$G$ can be modified into $H$ by a sequence of vertex deletions and local complementation. 
If $G$ has no vertex-minor isomorphic to $H$, then $G$ is {\it $H$-vertex-minor-free}.
Note that if~$G$ contains~$H$ as a pivot-minor, then~$G$ contains~$H$ as a vertex-minor (but not necessarily vice versa).

Bouchet~\cite{Bo94} characterized circle graphs in terms of forbidden vertex-minors. By using this result, Geelen and Oum~\cite{GO09} were able to characterize circle graphs in terms of forbidden pivot-minors.
Geelen, Kwon, McCarty, and Wollan~\cite{GKMW19} proved that for every circle graph~$H$, the class of $H$-vertex-minor-free graphs has bounded rank-width.
In contrast, it can be easily observed that for every non-bipartite graph $H$, the class of $H$-pivot-minor-free graphs has unbounded rank-width (see~\cite{DDJKKOP21}).
Oum~\cite{Ou09} conjectured that for each fixed {\it bipartite} circle graph~$H$, the class of $H$-pivot-minor-free graphs has bounded rank-width.
This conjecture is known to be true for bipartite graphs~\cite{Ou05}, line graphs~\cite{Ou09} and circle graphs~\cite{Ou09}.

Kant\'e and Kwon~\cite{KK18} conjectured that for every tree $T$, the class of $T$-vertex-minor-free graphs has bounded linear rank-width. This conjecture holds if $T$ is a path~\cite{KMOW21}, and
even if $T$ is a caterpillar, as we observed in~\cite{DDJKKOP21} (a caterpillar is a a tree that contains a path~$P$, such that every vertex not on~$P$ has a neighbour in~$P$).
Moreover, the conjecture has been confirmed for every class of graphs whose prime graphs (with respect to split decompositions) have bounded linear rank-width~\cite{KK18}. 

Recently, we showed that if~$T$ is a tree that is not a caterpillar, then the class of $T$-pivot-minor-free distance-hereditary graphs has unbounded linear rank-width~\cite{DDJKKOP21}. This led us to conjecture that for every caterpillar $T$, the class of $T$-pivot-minor-free graphs has bounded linear rank-width. In the same paper~\cite{DDJKKOP21}, we 
showed that this conjecture holds for the class of distance-hereditary graphs, and also confirmed the conjecture if $T$ has at most four vertices.
Hence, the class of claw-pivot-minor-free graphs has bounded linear rank-width, while its shrub-depth, equivalently rank-depth, is unbounded~\cite{KMOW21}.

Dahlberg, Helsen, and Wehner~\cite{Dahlberg2022} showed that if $G$ and $H$ are part of the input, then the problem of deciding if $G$ contains a vertex-minor isomorphic to $H$ is \NP-complete. 

\subsection{Our Results}

So far, results on pivot-minors were mainly of a structural nature. In our paper, we study pivot-minors from an \emph{algorithmic} perspective and perform, for the first time, a systematic study into the complexity of computing pivot-minors. That is, we consider the following research question:\\[0.5em]
{\em Can we decide in polynomial time whether a graph~$H$ is isomorphic to a pivot-minor of a graph~$G$?}\\[0.5em]
We denote the corresponding decision problem as:

\problemdef{\sc Pivot-Minor}{A pair of graphs~$G$ and~$H$.}{Does~$G$ contain a pivot-minor isomorphic to $H$?}

\noindent
In Section~\ref{s-npc} we prove that {\sc Pivot-Minor} is \NP-complete.
As a natural next step, we consider in the remainder of our paper the direction proposed in~\cite{Ou16} (see Question~7), that is, we assume that $H$ is not part of the input but fixed in advance. In this case, we denote the problem as {\sc $H$-Pivot-Minor}.

We first observe that for every graph~$H$, {\sc $H$-Pivot-Minor} is polynomial-time solvable for graphs of bounded rank-width. 
The reason for this result is that pivot-minor testing can be expressed in monadic second-order logic with modulo-$2$ counting in a similar way as done for vertex-minor testing~\cite{CO07}, by modifying the formula in \cite[Theorem 6.5]{CO07} to force $Z_a$ empty as described by the concept of $\alpha\beta$-minors of isotropic systems in \cite{Ou08}.
We recall that for every caterpillar $H$ on at most four vertices, the class of $H$-pivot-minor-free graphs has bounded linear rank-width, and thus bounded rank-width. 
Hence, for every such graph~$H$, we find that {\sc $H$-Pivot-Minor} is polynomial-time solvable.
 
Just like as for our other polynomial-time results, we show that for these cases we can even find a {\it certifying} polynomial-time algorithm. We discuss this in more detail below, but first
we summarize our results on {\sc $H$-Pivot-Minor} in the following state-of-the-art theorem. Here, we write $F\ssi G$ to denote that $F$ is an induced subgraph of $G$ (see Section~\ref{s-prelim} for other notation used below).

\begin{theorem}\label{t-main}
    Let $H$ be a graph not in $\{K_4,C_3+P_1\}$ such that 
    \begin{itemize}
    \item $H\ssi P_3+tP_1$ for some integer $t\geq 0$, or
    \item $H=K_{1,t}$ for some integer $t\geq 1$, or
    \item $|V(H)|\leq 4$.
\end{itemize}
Then there is a \emph{certifying} algorithm that solves $H$-{\sc Pivot-Minor} in polynomial time.
\end{theorem}

To explain the idea behind our algorithms, we  
observe that for every graph~$H$, the class of $H$-pivot-minor-free graphs is closed under vertex deletions.
Let~${\cal F}_H$ be the set of induced-subgraph-minimal graphs that contain a pivot-minor isomorphic to $H$.
Then, a graph~$G$ is $H$-pivot-minor-free if and only if it contains no induced subgraph in ${\cal F}_H$.

In Section~\ref{s-tp1free}, we show that if $H\ssi P_3+tP_1$ for some $t\geq 0$ or $H=K_{1,t}$ for some $t\geq 1$, then there exists an integer $c_H$ such that all graphs in ${\cal F}_H$ have at most $c_H$ vertices. 
In Section~\ref{s-small}, we show that 
if $H\notin \{C_3, \text{paw}, \text{diamond}, K_4,C_3+\nobreak P_1\}$ with $|V(H)|\leq 4$, then there exists an integer $c_H$ such that all graphs in ${\cal F}_H$ have at most $c_H$ vertices. 
Hence, for every such graph~$H$, we can solve {\sc $H$-Pivot-Minor} in polynomial time: it suffices to check by brute force if the input graph $G$ contains a graph from ${\cal F}_H$ as an induced subgraph after generating all non-isomorphic graphs in  $\mathcal{F}_H$. This algorithm is {\it certifying}: if it finds a graph $F$ from ${\cal F}_H$, then we can use $F$ as a certificate, as $F$ contains a pivot-minor isomorphic to $H$ by definition. 
See~\cite{MMNS11} for a survey on certifying algorithms.

In Section~\ref{s-small}, we also show that the class of $C_3$-pivot-minor-free graphs coincides with the class of bipartite graphs. Hence, as a certificate we can take any odd cycle, which we can readily find in polynomial time. Moreover, in the same section, we show that a graph is paw-pivot-minor-free if and only if it is diamond-pivot-minor-free. We prove that ${\cal F}_{paw}$ consists of the paw, the diamond, and all odd holes (cycles of length at least~$5$). We obtain a certifying algorithm by analyzing a structural property. Alternatively, one can use the more general algorithm of Chudnovsky, Scott, Seymour and Spirkl~\cite{CSSS20} for detecting an odd hole. 

There are two remaining graphs $H$ on at most four vertices: $H=K_4$ and $H=C_3+P_1$. We do not know the computational complexity of 
{\sc $H$-Pivot Minor} if $H\in \{K_4,C_3+P_1\}$. 
In Section~\ref{s-infi}, we prove that~${\cal F}_{K_4}$ and~${\cal F}_{C_3+\nobreak P_1}$ each contain infinitely many non-isomorphic graphs.

In Section~\ref{s-claw} we return to the case $H=K_{1,3}$. We showed that all graphs in ${\cal F}_{K_{1,t}}$ have bounded size for every $t\geq 1$, but in this result we do not specify the list of forbidden induced subgraphs.
However, for $t=3$, we are able to provide an explicit list, as we already showed in the conference version~\cite{Dabrowski2018} of our paper. Since this result has been 
explicitly mentioned in~\cite{DDJKKOP21,KMOW21}, we included its proof in Section~\ref{s-claw}. 
In Section~\ref{s-con} we conclude our paper with a brief discussion on future research directions.

\section{Preliminaries}\label{s-prelim}

All graphs $G=(V(G),E(G))$ have no loops and no multiple edges. We may write $V=V(G)$ and $E=E(G)$.
	For a subset $S\subseteq V(G)$, let~$G[S]=(S,\{uv\; |\; uv\in E, u,v\in S\})$ denote the subgraph of~$G$ {\it induced by}~$S$. 
	A graph~$H$ is an \emph{induced subgraph} of $G$ if $H=G[S]$ for some $S\subseteq V(G)$.
	For a vertex $v\in V(G)$, we let $G- v$ be the graph obtained from~$G$ by removing~$v$. 
	For a set $S\subseteq V(G)$, we let $G-S$ be the graph obtained from~$G$ by removing all vertices in~$S$. 
	For an edge $e\in E(G)$, we let $G-e$ be the graph obtained from~$G$ by removing~$e$. For  
	a set $F\subseteq E(G)$, we let $G-F$ be the graph obtained from~$G$ by removing all edges in $F$.
	
Let $G=(V,E)$. For two disjoint vertex subsets $A$ and $B$ in $G$, we say that $A$ is \emph{complete} to $B$ if every $a\in A$ is adjacent to every $b\in B$, while $A$ is \emph{anti-complete} to $B$ if every $a\in A$ is non-adjacent to every $b\in B$. A vertex~$u\in V$  is {\it (anti)-complete} to a set~$S\subseteq V\setminus \{u\}$ if $u$ is (non-)adjacent to every vertex of $S$. 
Two vertices $u$ and $v$ of $G$ are {\it twins} if for every $w\in V\setminus \{u,v\}$, it holds that $w$ is adjacent to $u$ if and only if $w$ is adjacent to $v$.
    
The graph~$\overline{G}=(V,\{uv\; |\; uv\notin E(G), u\neq v\}$ is the \emph{complement} of a graph $G$.
The graph $G_1+\nobreak G_2=(V(G_1)\cup V(G_2),E(G_1)\cup E(G_2))$ is the \emph{disjoint union} of two vertex-disjoint graphs~$G_1$ and $G_2$.
For a set of graphs $\{H_1, \ldots, H_p\}$, a graph~$G$ is {\em $(H_1,\ldots,H_p)$-free} if~$G$ does not contain an induced subgraph isomorphic to any graph in~$\{H_1, \ldots, H_p\}$.

\begin{figure}
   \centering
   \begin{center}
     $\quad\quad$
       \begin{tikzpicture}
        \tikzstyle{v}=[circle,draw,fill=black!50,inner sep=0pt,minimum width=4pt]
       \draw (0,0) node[v](v){};
       \draw (1,.5) node[v](v1){};
       \draw(1,-.5)node[v](v2){};
       \draw (-1,0) node[v](v3){};
       \draw(v3)--(v)--(v1)--(v2)--(v);
       \draw(0,-.5) node [label=below:paw]{};
     \end{tikzpicture}
      $\quad\quad$
          \begin{tikzpicture}
            \tikzstyle{v}=[circle,draw,fill=black!50,inner sep=0pt,minimum width=4pt]
       \draw (0,0) node[v](w2){};
       \draw (1,.5) node[v](v1){};
       \draw(1,-.5)node[v](v2){};
       \draw (2,0) node[v](w1){};
       \foreach \i in {1,2}{
         \draw (w1)--(v\i)--(w2);
       }
       \draw (v1)--(v2);
       \draw(1,-.5) node [label=below:diamond]{};
     \end{tikzpicture}
     $\quad\quad$
     \begin{tikzpicture}
      \tikzstyle{v}=[circle,draw,fill=black!50,inner sep=0pt,minimum width=4pt]
       \draw (0,0) node[v](v){};
       \draw (.5,.5) node[v](v1){};
       \draw(.5,-.5)node[v](v2){};
       \draw (-1,0) node[v](v3){};
       \draw(1,0)node[v](w){};
       \draw(v3)--(v)--(v1)--(w)--(v2)--(v)--(w);
       \draw(0,-.5) node [label=below:dart]{};
     \end{tikzpicture}
         $\quad\quad$
               \begin{tikzpicture}
                \tikzstyle{v}=[circle,draw,fill=black!50,inner sep=0pt,minimum width=4pt]
     \draw (0,0) node [v] (v1){};
     \draw (1,0) node [v] (v2){};
     \draw (1,1) node [v] (v3){};
     \draw (0,1) node [v] (v){};
     \foreach \i in {1,2,3}{
       \draw (v)--(v\i);
     }
     \draw (0.5,0) node[label=below:claw]{};
   \end{tikzpicture}
 \\[5pt]
    \begin{tikzpicture}
      \tikzstyle{v}=[circle,draw,fill=black!50,inner sep=0pt,minimum width=4pt]
    \draw (0,0) node [v] (v1){};
    \draw (1,0) node [v] (v2){};
    \draw (1,1) node [v] (v3){};
    \draw (0,1) node [v] (v4){};
    \draw (.5,.5) node[v] (v) {};
    \draw (v4)--(v1)--(v)--(v2)--(v3);
    \draw (v1)--(v2);
    \draw (0.5,0) node[label=below:bull]{};
  \end{tikzpicture}
  $\quad\quad$
  \begin{tikzpicture}
    \tikzstyle{v}=[circle,draw,fill=black!50,inner sep=0pt,minimum width=4pt]
    \draw (0,0) node [v] (v1){};
    \draw (1,0) node [v] (v2){};
    \draw (1,1) node [v] (v3){};
    \draw (0,1) node [v] (v4){};
    \draw (.5,.5) node[v] (v) {};
    \foreach \i in {1,2,3,4} {
      \draw (v)--(v\i);
      }
      \draw (v1)--(v2)--(v3)--(v4)--(v1);
      \draw (0.5,0) node[label=below:$W_4$]{};
  \end{tikzpicture}
  $\quad\quad$
    \begin{tikzpicture}
      \tikzstyle{v}=[circle,draw,fill=black!50,inner sep=0pt,minimum width=4pt]
    \begin{scope}[xshift=-1.5cm]
      \foreach \i in {1,2} {
      \draw (120*\i:.7) node[v,label=left:$a_\i$](a\i){};
    }
      \foreach \i in {3} {
      \draw (-10:.7) node[v,label=left:$a_\i$](a\i){};
    }
    \draw (a1)--(a2)--(a3)--(a1);
    \end{scope}
    \begin{scope}[xshift=1.5cm]
      \draw (10:1) node[v,label=right:$b_4$] (b4){};
      \foreach \i in {1,2} {
      \draw (180-120*\i:.7) node[v,label=right:$b_\i$](b\i){};
      \draw (a\i)--(b\i)--(b4);
    }
      \foreach \i in {3} {
      \draw (190:.7) node[v,label=above left:$b_\i$](b\i){};
      \draw (a\i)--(b\i)--(b4);
    }
    \draw (b1)--(b2)--(b3)--(b1);
    \end{scope}
    \draw (-90:.4) node[label=below:$\bw_3$]{};
  \end{tikzpicture}
    \begin{tikzpicture}
      \tikzstyle{v}=[circle,draw,fill=black!50,inner sep=0pt,minimum width=4pt]
    \draw (0,0) node [v] (a1){};
    \draw (0.5,0.5) node [v] (a2){};
    \draw (0,1) node [v] (a3){};

    \draw (2,0) node [v] (b1){};
    \draw (1.5,0.5) node [v] (b2){};
    \draw (2,1) node [v] (b3){};

    \draw(a1)--(a2)--(a3)--(a1);
    \draw(b1)--(b2)--(b3)--(b1);
    \draw(a1)--(b1);
    \draw(a2)--(b2);
    \draw(a3)--(b3);
  
\draw (1,0) node[label=below:prism]{};
  \end{tikzpicture}
   \end{center}
 \caption{The paw, diamond, dart, claw, bull, $W_4$, $\overline{BW_3}$, and prism.}
 \label{fig:s2}
 \end{figure}
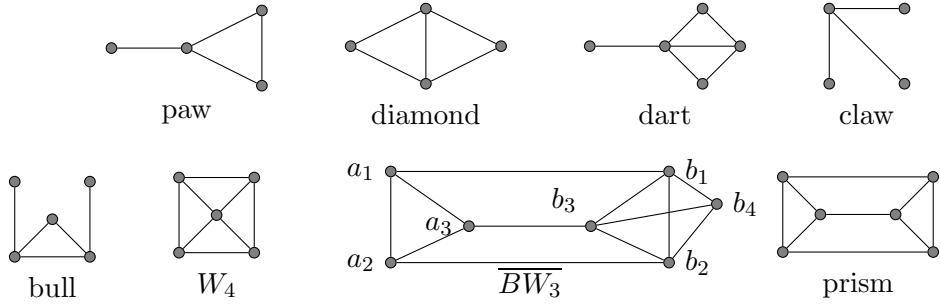

The path, cycle, star and complete graph on~$n$ vertices are denoted by~$P_n$, $C_n$, $K_{1,n-1}$ and~$K_n$, respectively.
The graph~$W_n$ is obtained from~$C_n$ by adding one vertex adjacent to all vertices of~$C_n$.
The \emph{paw}, \emph{diamond}, \emph{dart} and \emph{claw} are the graphs $\overline{P_1+P_3}$, $\overline{2P_1+P_2}$, $\overline{P_1+\text{paw}}$ and~$K_{1,3}$, respectively.
The \emph{bull} is the graph obtained from~$P_5$ by adding an edge between the second vertex and the fourth vertex.
The \emph{prism} is the complement of $C_6$.
The graph~$BW_3$ is the bipartite graph on seven vertices obtained from~$C_6$ by adding one vertex adjacent to three pairwise non-adjacent vertices of the cycle. In our proof, the graph $\overline{BW_3}$ will play an important role. 
See Figure~\ref{fig:s2} for an illustration of the above graphs.
 
Two graphs are \emph{pivot-equivalent} if one can be obtained from the other by a sequence of edge pivots. A graph class is \emph{pivot-minor-closed} if it is closed under vertex deletions and edge pivots.
A sequence~$S$ of vertex deletions and edge pivots is an \emph{$H$-pivot-minor-sequence} of a graph~$G$ if~$H$ can be obtained from~$G$ after applying the operations of~$S$.

We write~$G/v$ to denote $(G\pivot zv)- v$ if a vertex~$v$ has a neighbour~$z$ and~$G-v$ if~$v$ is isolated. 
This is a well-defined notation up to the pivot-equivalence. Namely,
for two distinct neighbours~$x$,~$y$ of~$v$, it holds that $(G\pivot xv)- v$ is pivot-equivalent to $(G\pivot yv)- v$, as $$(G\pivot xv)-v=(G\pivot yv\pivot xy)-v=(G\pivot yv-v)\pivot xy.$$ Hence, the choice 
of neighbour of $v$ does not change the pivot-equivalence of graphs~$G/v$.

We need two basic lemmas. Lemma~\ref{lem:bouchet} holds in the context of binary delta-matroids or matrix pivots (see~\cite{BD91,Ou09}).
We provide a direct proof, inspired by the analogous proof for vertex-minors~in~\cite{GO09}.

\begin{lemma}\label{lem:backward}
Let $v,x,y$ be distinct vertices of a graph~$G$.
If~$xy\in E(G)$, then
$(G\pivot xy)- v$ is pivot-equivalent to~$G-v$ and $(G\pivot xy)/v$ is pivot-equivalent~to~$G/v$.
\end{lemma}
\begin{proof}
The first statement is trivial.
We prove the second statement.
As a pivot operation does not change the fact that~$v$ has a neighbour in~$G$, $v$ has a neighbour in~$G$ if and only if it has a neighbour in~$G\pivot xy$.
If~$v$ is isolated in~$G$, then the result follows from the first statement.
We may assume~$v$ has a neighbour in~$G$ and in~$G\pivot xy$.
If~$v$ is non-adjacent to~$x$ and~$y$ in~$G$, then for a neighbour~$z$ of~$v$, we have that $(G\pivot xy)/v=(G\pivot xy\pivot zv)- v= (G\pivot zv- v)\pivot xy=(G/v)\pivot xy$.
Thus~$(G\pivot xy)/v$ is pivot-equivalent to~$G/v$.
If~$v$ is adjacent to~$x$, then $(G\pivot xy) /v $ is pivot-equivalent to $(G\pivot xy \pivot yv)-v = (G\pivot xv)-v=G/v$.
\end{proof}

\begin{lemma}\label{lem:bouchet}
If a graph~$H$ is a pivot-minor of a graph $G$ and $v\in V(G)\setminus V(H)$, then~$H$ is a pivot-minor of~$G-v$ or $(G\pivot vw)- v$ for some neighbour~$w$ of~$v$~in~$G$.
\end{lemma}
\begin{proof}
Let $G_0=G$, $G_1=G_0\pivot x_1y_1$, $G_2=G_1\pivot x_2y_2$, $\ldots$, $G_m=G_{m-1}\pivot x_my_m$ and suppose~$H$ is an induced subgraph of~$G_m-v$.
For each~$i$, if $x_i\neq v$ and $y_i\neq v$ then by Lemma~\ref{lem:backward}, $G_i-v=(G_{i-1}\pivot x_iy_i)-v$ is pivot-equivalent to~$G_{i-1}-v$ and $G_i/v=(G_{i-1}\pivot x_iy_i)/v$ is pivot-equivalent to~$G_{i-1}/v$.
If $x_i=v$, then $G_i-v=(G_{i-1}\pivot x_iy_i)-v =G_{i-1}/v$ and $G_i/v=(G_{i-1}\pivot y_iv \pivot y_i v) -v = G_{i-1}-v$.
Thus we deduce that~$G/v$ or~$G-v$ is pivot-equivalent to~$G_m-v$.
\end{proof}

\section{When~$\mathbf{H}$ Is Part of the Input}\label{s-npc}

In this section, we prove that {\sc Pivot-Minor} is \NP-complete.
We first introduce some terminology and basic results on matroids.
A \emph{matroid} is a pair $M=(E,{\cal I})$ of a finite set~$E$, called the \emph{ground set}, and a set~$\mathcal I$ of subsets of~$E$ satisfying 
the following three properties: 
\begin{itemize}
\item ${\cal I}\neq \emptyset$;
\item if $Y\in {\cal I}$ and $X\subseteq Y$, then $X\in {\cal I}$ and
\item if $X,Y\in {\cal I}$ with $|Y|=|X|+1$, then there exists an element $y\in Y\setminus X$ such that $X\cup \{y\}\in {\cal I}$.
\end{itemize}
A set $X\subseteq E$ is \emph{independent} in $M=(E,\mathcal I)$ if $X\in \mathcal I$.
Otherwise~$X$ is \emph{dependent}.
The \emph{rank} of a subset $X\subseteq E$ is the size of a largest independent subset of~$X$.
The \emph{rank} of a matroid $M=(E,\mathcal I)$ is the rank of~$E$.
A \emph{base} of a matroid is a maximal independent set.
A \emph{circuit} of a matroid is a minimal dependent set.
The \emph{dual matroid}~$M^*$ of a matroid $M=(E,\mathcal I)$ is a matroid on~$E$ such that~$X$ is a base of~$M^*$ if and only if~$E\setminus X$ is a base in~$M$.
For a subset~$X$ of~$E$, we define $M\setminus X$ to be the matroid $(E\setminus X, \mathcal I')$ such that $\mathcal I'=\{ X'\subseteq E\setminus X\; |\; X'\in \mathcal I\}$.
We define $M/X=(M^*\setminus X)^*$.
A matroid~$N$ is a \emph{minor} of a matroid~$M$
if $N=(M\setminus X)/Y$
for some disjoint sets~$X$ and~$Y$.
A matroid $M=(E,\mathcal I)$ is \emph{binary} if there is a matrix over the binary field whose columns are indexed by~$E$ such that~$X$ is independent in~$M$ if and only if the corresponding columns are linearly independent.
It is known that the dual matroid of a binary matroid is also binary.

A major example of binary matroids arises from graphs.
For a graph $G=(V,E)$, let~$\mathcal I$ be the set of subsets~$X$ of~$E$ such that the subgraph $(V,X)$ has no cycles.
Then $M(G)=(E,\mathcal I)$ is a matroid, called the \emph{cycle matroid} of~$G$ and such matroids are binary.
It is known that circuits of~$M(G)$ are precisely the edge set of cycles of~$G$. We also use the following fact.

\begin{lemma}[see~\cite{Oxley}]\label{lem:minor}
    If a graph~$H$ is a minor of a graph~$G$, then~$M(H)$ is a minor of~$M(G)$.
\end{lemma}

If~$G$ is connected and has~$n$ vertices and~$m$ edges, then~$M(G)$ has rank $n-1$ because any spanning tree of~$G$ has $n-1$ edges, and~$(M(G))^*$ has rank $m-n+1$.

For a binary matroid $M=(E,\mathcal I)$, the \emph{fundamental graph} of~$M$ with respect to a base~$B$ is the bipartite graph on~$E$ with the bipartition $(B, E\setminus B)$ such that $x\in B$, $y\in E\setminus B$ are adjacent if and only if $(B\setminus\{x\})\cup \{y\}$ is a base of~$M$.
Conversely, for a bipartite graph~$G$ with a bipartition~$(A,B)$, we may define a binary matroid $\Bin(G,A,B)$ on~$V(G)$ represented by the $A\times V(G)$ matrix
\[
\bordermatrix{
& A & B \cr
A & I_A & M_{A,B}
}
\] over the binary field where~$I_A$ is the $A\times A$ identity matrix and~$M_{A,B}$ is the $A\times B$ submatrix of the adjacency matrix of~$G$ whose $(x,y)$-entry is~$1$ if and only if~$x$ and~$y$ are adjacent.
We need the following lemma for our \NP-hardness result.

\begin{lemma}[{\cite[Corollary~3.6]{Ou05}}]\label{lem:bin}
The following statements hold:
\begin{enumerate}[label=\rm(\roman*)]
\item Let~$N$ and $M$ be binary matroids, and let~$H$ and $G$ be fundamental graphs of~$N$ and~$M$ respectively.
If~$N$ is a minor of~$M$, then~$H$ is a pivot-minor of~$G$.
\item Let~$G$ be a bipartite graph with bipartition $(A, B)$.
If~$H$ is a pivot-minor of~$G$, then there is a bipartition $(A', B')$ of $H$ such that $\Bin(H,A',B')$ is a minor of $\Bin(G,A,B)$.
\end{enumerate}
\end{lemma}

We are now ready to prove our hardness result.

\begin{theorem}\label{t-hard}
{\sc Pivot-Minor} is \NP-complete.
\end{theorem}

\begin{proof}
We reduce from the {\sc Hamiltonian Cycle} problem, which asks if a  graph has a Hamiltonian cycle.
This problem is \NP-complete even for $3$-regular graphs~\cite{GJT76}.
Let $G=(V,E)$ be a $3$-regular graph with~$n$ vertices and~$m$ edges.
We may assume without loss of generality that $n\geq 5$ and that~$G$ is connected.
As $G$ is 3-regular, $2m=3n$.
Consequently, $(M(G))^*$ has rank $m - n + 1 =\frac{1}{2}n+1$.

Let~$T$ be a spanning tree of~$G$. Let
$G_T$ be the fundamental graph of~$M(G)$ with respect to~$E(T)$, which can be built in polynomial time.
We claim that
\begin{itemize}
    \item $G$ has a Hamiltonian cycle if and only if $G_T$ contains a pivot-minor isomorphic to~$K_{1,n-1}$.
\end{itemize}

For the forward direction, we use Lemma~\ref{lem:bin}(i).
If~$G$ has a Hamiltonian cycle~$C$, then~$G$ contains~$C$ as a minor and by Lemma~\ref{lem:minor}, $M(G)$ has~$M(C)$ as a minor. So~$G_T$ has every fundamental graph of~$M(C)$ as a pivot-minor.
This proves the forward direction because every fundamental graph of~$M(C)$ is isomorphic to~$K_{1,n-1}$.

For the reverse direction, suppose that~$G_T$ contains a pivot-minor isomorphic to~$K_{1,n-1}$.
Then by Lemma~\ref{lem:bin}(ii), $V(K_{1,n-1})$ has a bipartition $(A',B')$ such that $\Bin(K_{1,n-1},A',B')$ is a minor of
$M(G)=\Bin(G_T,A,B)$ for some partition~$(A,B)$ of~$V(G_T)$.

As~$K_{1,n-1}$ is connected, the center of $K_{1,n-1}$ forms one part of $(A',B')$ and the rest forms the other part. 
So $\Bin(K_{1,n-1},A',B')$ is either~$M(C)$ or its dual~$(M(C))^*$, where~$C$ is the cycle on~$n$ vertices.
Therefore~$M(C)$ or~$(M(C))^*$ is a minor of~$M(G)$.
Equivalently, $M(C)$ is a minor of~$M(G)$ or~$(M(G))^*$.
Because the rank of~$M(C)$ is $n-1$ and the rank of~$(M(G))^*$ is $\frac12 n+1 < n-1$ (as $n\geq 5$) we find that~$M(C)$ cannot be a minor of~$(M(G))^*$.
Thus, $M(C)$ is a minor of~$M(G)$ and therefore~$M(G)$ has a circuit of length at least~$n$.
This implies that~$G$ has a cycle of length~$n$.
\end{proof}

\section{When $\mathbf{H}$ Is Fixed: $\mathbf{H=tP_1, P_2+tP_1, P_3+tP_1}$, or $\mathbf{K_{1,t}}$}\label{s-tp1free}

In this section, we provide a certifying algorithm that solves \textsc{$H$-Pivot-Minor} in polynomial time when $H=tP_1$, $H=P_2+tP_1$, $H=P_3+tP_1$, or $H=K_{1,t}$, for every $t\geq 1$.
As explained in Section~\ref{s-intro}, for every such graph $H$, our approach is to show that there exists an integer~$c_H$ such that all graphs in  ${\cal F}_H$ have at most $c_H$ vertices.
To prove this, we adjust a concept for vertex-minors introduced by Geelen and Oum~\cite{GO09} to pivot-minors.
Namely, a graph $G$ is \emph{$H$-pivot-unique} if $G$ contains a pivot-minor isomorphic to $H$ and  for every vertex $v$ of $G$, at most one of $G-v$ and $G/v$ contains a pivot-minor isomorphic to $H$.
The following lemma is readily seen.

\begin{lemma}\label{lem:pivotunique}
   For every graph~$H$, if a graph $G$ contains a pivot-minor isomorphic to $H$ and every proper induced subgraph of $G$ is $H$-pivot-minor-free,
    then $G$ is $H$-pivot-unique.
\end{lemma}

We first prove a sequence of five general lemmas, starting with the following lemma, which states that $H$-pivot-uniqueness is preserved under pivoting.

\begin{lemma}\label{lem:pivot-pu}
For every graph~$H$, if a graph $G$ is $H$-pivot-unique,  then every graph that is pivot-equivalent to~$G$ is $H$-pivot-unique.
\end{lemma}
\begin{proof}
    Assume a graph $G$ is $H$-pivot-unique. Let $v\in V(G)$ and $wz\in E(G)$. It suffices to show that if at most one of $G-v$ and $G/v$ contains a pivot-minor isomorphic to $H$, then at most one of $(G\pivot wz)-v$ and $(G\pivot wz)/v$ contains a pivot-minor isomorphic to $H$.
    Suppose that at most one of $G-v$ and $G/v$ contains a pivot-minor isomorphic to $H$.
    If $v\in \{w,z\}$, then this holds because $G-v$ is pivot-equivalent to $(G\pivot wz)/v=(G\pivot wz\pivot wz)-v$, and 
    $G/v=(G\pivot wz)-v$ is pivot-equivalent to $(G\pivot wz)-v$.
    If $v\notin \{w,z\}$, then we apply Lemma~\ref{lem:backward}.
\end{proof}

\begin{lemma}\label{lem:pivot-twin}
    Let $H$ be an induced subgraph of a graph $G$. 
    Let $u$ and $v$ be adjacent vertices of~$G$ with $v\notin V(H)$.
    If $u$ and $v$ are twins in $G[V(H)\cup\{u,v\}]$,
    then $G$ is not $H$-pivot-unique.
\end{lemma}
\begin{proof}
    As $G[V(H)\cup\{u,v\}]\pivot uv=G[V(H)\cup\{u,v\}]$, 
    both $G/v$ and $G-v$ contain $H$ as a pivot-minor.
\end{proof}

\begin{lemma}\label{lem:truetwin}
    Let $G$ be a graph, $xy\in E(G)$, and $v,w$ be adjacent twins in $G$. Then $v$ and $w$ are adjacent twins in $G\pivot xy$. 
\end{lemma}
\begin{proof}
    If $\{x,y\}=\{v,w\}$ or $\{x,y\}\cap \{v,w\}=\emptyset$, then it is obvious. We assume that $\abs{\{x,y\}\cap \{v,w\}}= 1$. Without loss of generality, assume that $w=x$.

    Note that $v$ is a common neighbour of $x$ and $y$. The set $(N_G(x)\setminus N_G(y))\setminus \{y\}$ is complete to~$x$ in $G$, and thus it is anti-complete to $x$ in $G\pivot xy$.
    The set $(N_G(y)\setminus N_G(x))\setminus \{x\}$ is anti-complete to $x$ in $G$, and thus it is complete to $x$ in $G\pivot xy$. Therefore, $v$ is adjacent to $w=x$ in $G\pivot xy$, and $N_{G\pivot xy}(v)\setminus \{w\}=N_{G\pivot xy}(w)\setminus \{v\}=(N_G(y)\setminus \{v,x\})\cup \{y\}$. This shows the lemma.
\end{proof}

\begin{lemma}\label{lem:pivot-isolated}
    Let $H$ be an induced subgraph of a graph $G$. 
    If there is a vertex~$v$ in $V(G)\setminus V(H)$ with no neighbours in $H$,
    then $G$ is not $H$-pivot-unique.
\end{lemma}
\begin{proof}
    This is trivial because $G/v$ contains $H$ as an induced subgraph.
\end{proof}

We now consider the case where $H=tP_1$ for some $t\geq 1$.

\begin{theorem}\label{thm:tP1}
   For every $t\geq 1$, if a graph $G$ is $tP_1$-pivot-unique, 
    then $\abs{V(G)}\le 2^{t}-1$.
\end{theorem}
\begin{proof}
    Let $G$ be a $tP_1$-pivot-unique graph.
    By Lemma~\ref{lem:pivot-pu}, we may assume that 
    $tP_1$ is an induced subgraph of $G$.
    Let $X$ be an independent set of $t$ vertices in $G$.
    For a subset $Y$ of $X$, let $N_Y$ be the set of vertices~$v$ in $V(G)\setminus X$ such that $N_G(v)\cap X=Y$.

    By Lemma~\ref{lem:pivot-twin}, 
    $N_Y$ is independent for all $Y\subseteq X$.
    Furthermore, again by Lemma~\ref{lem:pivot-twin},
    $N_Y=\emptyset$ if $\abs{Y}=1$.
    By Lemma~\ref{lem:pivot-isolated}, $N_\emptyset=\emptyset$.

    We claim that  $\abs{N_Y}\le 1$ for each non-empty $Y\subseteq X$. 
    Let $Y$ be a non-empty subset of $X$, 
    $u$ be a vertex in $Y$, and 
    $w\in N_Y$.
    Then in $G\pivot uw-w$, $(N_Y\cup X)\setminus \{w,u\}$ is independent. 
    Since~$G$ is $tP_1$-pivot-unique, 
    $\abs{(N_Y\cup X)\setminus \{w,u\}}<t$.
    This proves that $\abs{N_Y}\le 1$ for each non-empty~$Y\subseteq X$.
    This implies that $\abs{V(G)}\le t + (2^{t}-1-t) = 2^t -1$.
\end{proof}

We further show that 
    if $G$ is an induced-subgraph-minimal graph having a pivot-minor isomorphic to $P_2+tP_1$, then 
    $\abs{V(G)}\le (t+1)(2^{t+2}-t-2)$. In this case, the property of being $(P_2+tP_1)$-pivot-unique does not directly bound the size of the graph, as $K_{1,n}+tP_1$ is $(P_2+tP_1)$-pivot-unique for every $n$. But $K_{1,n}+tP_1$ with large $n$ is not an induced-subgraph-minimal graph having a pivot-minor isomorphic to $P_2+tP_1$. So, we refine the argument.
 \begin{lemma}\label{lem:stariso1}
 Let $G$ be a $(P_2+tP_1)$-pivot-unique graph and $v$ be a vertex of~$G$.
  If $v$ has more than $t$ neighbours of degree~$1$, then
  $G$ is isomorphic to $K_{1,n} + mP_1$ for some integers $n$ and $m$.
 \end{lemma}
 
 \begin{proof}
  Let $A$ be the set of neighbours of~$v$ having degree~$1$.
   Let $w\in A$. 
   If $G-(A\cup \{v\})$ contains an edge, then 
   both $G-w$ and $(G\pivot vw)-w$ contain induced subgraphs isomorphic to $P_2+tP_1$, 
   contradicting the assumption that $G$ is $(P_2+tP_1)$-pivot-unique.
    So, $G-(A\cup \{v\})$ contains no edges. This implies that
    every edge of $G$ is incident with $v$ and so $G$ is isomorphic to $K_{1,n}+mP_1$ for some integers $n$ and $m$.
 \end{proof}

 \begin{lemma}\label{lem:stariso2}
    Let $G$ be a graph such that every connected component is a star or a complete graph. 
    Let $K$ be a graph.
    If $G$ has a pivot-minor isomorphic to $K$ but no proper induced subgraph of $G$ has a pivot-minor isomorphic to $K$, then $G$ is isomorphic to $K$.
 \end{lemma}
 \begin{proof}
    Observe that 
    if every connected component of a graph is a star or a complete graph, 
    then 
    pivoting any edge produces an isomorphic graph.
    Therefore, for every induced subgraph~$H$ of~$G$ and a graph~$K$, 
    $H$ has a pivot-minor isomorphic to $K$
    if and only if $H$ has an induced subgraph isomorphic to~$K$. 
    It follows that $G$ is isomorphic to $K$.
 \end{proof}

Now, we prove the following.

\begin{theorem}\label{thm:P2tP1}
    Let $t$ be a positive integer. If $G$ is an induced-subgraph-minimal graph having a pivot-minor isomorphic to $P_2+tP_1$, 
    then $\abs{V(G)}\le (t+1)(2^{t+2}-t-2)$.
\end{theorem}
\begin{proof}
    By Lemma~\ref{lem:pivotunique}, $G$ is $(P_2+tP_1)$-pivot-unique.
    Since $G$ has a pivot-minor isomorphic to $P_2+tP_1$, 
    there exists a graph~$H$ pivot-equivalent to $G$ such that $H$ contains an induced subgraph isomorphic to $P_2+tP_1$. 
    By Lemma~\ref{lem:pivot-pu}, $H$
    is $(P_2+tP_1)$-pivot-unique.
    Let $X\subseteq V(H)$ such that $H[X]$ is isomorphic to $P_2+tP_1$.
    For a subset $Y$ of $X$, let $N_Y$ be the set of vertices~$v$ in $V(H)\setminus X$ such that $N_H(v)\cap X=Y$.

    By Lemma~\ref{lem:pivot-twin}, 
    $N_Y$ is independent for all $Y\subseteq X$.
    By Lemma~\ref{lem:pivot-isolated}, $N_\emptyset=\emptyset$.

    Let $v_1$, $v_2$ be the ends of the unique edge of $H[X]$.
    We claim that 
    \begin{itemize}
        \item[($\ast$)] for every $Y\subseteq X$, if $Y\setminus \{v_1,v_2\}\neq\emptyset$, then $\abs{N_Y}\le \min(\abs{Y}-1,1)$.    
    \end{itemize}
    Assume that $Y\subseteq X$ contains an isolated vertex of $H[X]$.
    If $\abs{Y}=1$, then by Lemma~\ref{lem:pivot-twin}, $N_Y=\emptyset$, proving the claim. 
    So we may assume that $\abs{Y}>1$ and $\abs{N_Y}>1$.
    Let $y$ be an isolated vertex of $H[X]$ in $Y$ and let $x\in N_Y$.
    Then in $H\pivot xy$, there are no edges between $X\setminus \{y\}$ and $N_Y\setminus \{x\}$, and so $H[(X\cup N_Y)\setminus \{x,y\}]$  contains an induced subgraph isomorphic to $P_2+tP_1$. Thus, both $H-x$ and $(H\pivot xy)-x$ contain pivot-minors isomorphic to $P_2+tP_1$, contradicting the fact that $H$ is $(P_2+tP_1)$-pivot-unique.
    Thus, the claim holds.

    By Lemma~\ref{lem:pivot-twin}, 
    $N_{\{v_1,v_2\}}=\emptyset$.
    If both $N_{\{v_1\}}$ and $N_{\{v_2\}}$ are non-empty, then both $H-v_1$ and $H\pivot v_1v_2-v_1$ contain pivot-minors isomorphic to $P_2+tP_1$, contradicting the fact that $H$ is $(P_2+tP_1)$-pivot-unique. 
    Thus we conclude that 
    \begin{itemize}
        \item[($\ast\ast$)] $N_{\{v_1,v_2\}}$ is empty and at least one of $N_{\{v_1\}}$ and $N_{\{v_2\}}$ is empty. 
    \end{itemize}
    By symmetry, we may assume that $N_{\{v_2\}}=\emptyset$.
    By $(\ast)$ and $(\ast\ast)$, we have 
    \[
      \left\lvert\bigcup_{Y\subseteq X, Y\neq\{v_1\}} N_Y\right\rvert
    \le 
    2^{t+2}-4-t.
    \]

    If $v_1$ has more than $t$ neighbours of degree~$1$ in $H$, 
    then by Lemma~\ref{lem:stariso1}, $H$ is isomorphic to $K_{1,n}+mP_1$ for some integers $n$ and $m$. 
    Since $G$ is pivot-equivalent to $H$, $G$ is also isomorphic to $K_{1,n}+mP_1$, as an edge pivot to a star yields an isomorphic graph.  
    Then by Lemma~\ref{lem:stariso2}, $G$ is isomorphic to $P_2+tP_1$. Thus, we are done. 
    Therefore, we may assume that $v_1$ has at most~$t$ neighbours of degree~$1$ in $H$. 

    If
    \(\abs{N_{\{v_1\}}}-t
    > t(2^{t+2}-4-t)\),
    then 
    by the pigeonhole principle, 
    $H$ has a vertex~$q\in \bigcup_{Y\subseteq X,Y\neq \{v_1\}} N_Y $ such that $q$ has more than $t$ neighbors in $N_{\{v_1\}}$. As $N_{\{v_2\}}$ and $N_{\{v_1, v_2\}}$ are empty by $(\ast\ast)$, $q$ has a neighbour in $X\setminus \{v_1, v_2\}$.
    Let $v_3\in X\setminus\{v_1,v_2\}$ be a neighbour of~$q$ in $H$.
    Let $U=N_{\{v_1\}}\cap N_H(q)$
    and $z\in U$.
    Observe that in $H\pivot qz$, 
    $v_1$ is 
    adjacent to $v_3$ 
    and there are no edges between $\{v_1,v_3\}$
    and $U\setminus\{z\}$.
    Thus, both $H-q$ and $(H\pivot qz)-q$ contain induced subgraphs isomorphic to $P_2+tP_1$, contradicting the fact that $H$ is $(P_2+2P_1)$-pivot-unique.
  
    Therefore $\abs{N_{\{v_1\}}}
    \le t(2^{t+2}-4-t)+t$.
    Then $\abs{V(G)}\le (t+2)+(2^{t+2}-4-t)+t(2^{t+2}-4-t)+t
    = (t+1)(2^{t+2}-t-2)$.
  \end{proof}

\begin{lemma}\label{lem:largestar1}
    Let $G$ be a $K_{1,t}$-pivot-unique graph. 
    If $G$ has a vertex having at least $t$ neighbours of degree $1$, then $G$ is a star.
 \end{lemma}
 \begin{proof}
    Let $v$ be a vertex having at least $t$ neighbours of degree~$1$
    and let $A$ be the set of neighbours of $v$ having degree~$1$.
    
    If $G-(A\cup \{v\})$ contains an edge $uz$, then 
    both $G-u$ and $(G\pivot uz)-u$ contains an induced subgraph isomorphic to $K_{1,t}$ on $A\cup \{v\}$, 
    contradicting the assumption that $G$ is $K_{1,t}$-pivot-unique.
    So, $G-(A\cup \{v\})$ contains no edges. 
    Therefore every edge is incident with~$v$.
    Since $G$ is $K_{1,t}$-pivot-unique, it has no isolated vertex. 
    Thus, $G$ is isomorphic to $K_{1,n}$ for some integer $n\ge t$. 
 \end{proof}

\begin{theorem}\label{thm:star}
    Let $t\ge 2$ be an integer. If $G$ is an induced-subgraph-minimal graph having a pivot-minor isomorphic to $K_{1,t}$, 
    then $\abs{V(G)}
    \le(t^2-1)(2^{t+1}-t-3)+2t+2$.
\end{theorem}
\begin{proof}
    By Lemma~\ref{lem:pivotunique}, $G$ is $K_{1,t}$-pivot-unique.
    Let $H$ be a graph pivot-equivalent to $G$ such that $H[X]$ is isomorphic to $K_{1,t}$ for a set $X$ of vertices.
    By Lemma~\ref{lem:pivot-pu}, $H$ is $K_{1,t}$-pivot-unique.
    Let $w$ be the center of the star $H[X]$.
    For a subset $Y$ of $X$, let $N_Y$ be the set of vertices~$v$ in $V(H)\setminus X$ such that $N_H(v)\cap X=Y$.

    By Lemma~\ref{lem:pivot-twin}, 
    $N_Y$ is independent for all $Y\subseteq X$.
    By Lemma~\ref{lem:pivot-isolated}, $N_\emptyset=\emptyset$.

    We claim that 
    \begin{itemize}
        \item[($\ast$)] for all $Y\subseteq X$, if $Y\setminus \{w\}\neq\emptyset$, then $\abs{N_Y}\le t-1$.
    \end{itemize}
    Suppose that $Y\subseteq X$, $Y\setminus \{w\}\neq\emptyset$, and $\abs{N_Y}\ge t$.
    Let $z\in Y\setminus \{w\}$.
    
   If $w\in Y$, then in both $H-w$ and $(H\pivot wz)-w$, $N_Y\cup \{z\}$ induces a subgraph isomorphic to~$K_{1,\abs{N_Y}}$. Thus, we may assume that $w\notin Y$.

    If there is $u\in X\setminus(Y\cup\{w\})$, then in both $H-w$ and $(H\pivot uw)-w$, $N_Y\cup \{z\}$ induces a subgraph isomorphic to $K_{1,\abs{N_Y}}$.
    So, we may assume that $Y=X\setminus\{w\}$.

    Since $t\ge 2$, there is 
    a vertex $p\in Y\setminus\{z\}$.
    In $H-p$, $N_Y\cup \{z\}$ induces a subgraph isomorphic to $K_{1,\abs{N_Y}}$.
    In $(H\pivot wp)-p$, $N_Y\cup \{w\}$ 
    induces a subgraph isomorphic to $K_{1,\abs{N_Y}}$.
    This contradicts the assumption that $H$ is $K_{1,t}$-pivot-unique.
    This proves $(\ast)$.

    By Lemma~\ref{lem:pivot-twin}, $N_Y=\emptyset$ for all subsets $Y$ of $X$ with $w\in Y$ and $\abs{Y}=2$.
    Also by Lemma~\ref{lem:pivot-twin}, 
    $N_X=\emptyset$.
    By $(\ast)$, we have \[\left\lvert\bigcup_{Y\subseteq X, Y\neq \{w\}}N_Y\right\rvert\le (t-1)(2^{t+1}-t-3).\] 

    If $N_{\{w\}}$ contains at least $t$ vertices of degree~$1$ in $H$, 
    then by Lemma~\ref{lem:largestar1}, $H$ is isomorphic to $K_{1,n}$ for some integer~$n$. 
    Since $G$ is pivot-equivalent to $H$, $G$ is also isomorphic to $K_{1,n}$, as an edge pivot to a star yields an isomorphic graph.  
    Then by Lemma~\ref{lem:stariso2}, $G$ is isomorphic to $K_{1,t}$. 
    Therefore, we may assume that $N_{\{w\}}$ contains fewer than $t$ neighbours of degree~$1$ in~$H$.

    If
    $\abs{N_{\{w\}}}-(t-1)
    > t(t-1)(2^{t+1}-t-3)$,
    then 
    by the pigeonhole principle, 
    there exists a vertex~$q\in\bigcup_{Y\subseteq X,Y\neq\{w\}}N_Y$ such that 
    $\abs{N_H(q)\cap N_{\{w\}}}>t$.
    Let $U=N_H(q)\cap N_{\{w\}}$ and $r\in U$.
    Observe that in $(H\pivot qr)-r$, 
    $(U\setminus\{r\})\cup\{w\}$ induces a subgraph isomorphic to $K_{1,{\abs{U}-1}}$. Thus, both $H-r$ and $(H\pivot qr)-r$ contain pivot-minors isomorphic to $K_{1,t}$, contradiciting the fact that $H$ is $K_{1,t}$-pivot-unique.

    Therefore $\abs{N_{\{w\}}}-(t-1)
    \le t(t-1)(2^{t+1}-t-3)$.
    Then $\abs{V(G)}
    \le (t+1)+(t+1)+(t-1)(2^{t+1}-t-3)+t(t-1)(2^{t+1}-t-3)=
    (t^2-1)(2^{t+1}-t-3)+2(t+1)$.
\end{proof}

\begin{lemma}\label{lem:P3iso1}
 Let $G$ be a $(P_3+tP_1)$-pivot-unique graph and $v$ be a vertex of~$G$.
  If $v$ has more than $t$ neighbours of degree~$1$, then
  $G$ is isomorphic to a graph obtained from $K_{1,n} + mP_1$ for some integers $n\ge t+1$ and $m\ge 0$, by 1-subdividing at most $n-(t+1)$ edges.
 \end{lemma}
 \begin{proof}
    Let $A$ be the set of neighbours of $v$ having degree~$1$. Let $w\in A$.
    
    If $G-(A\cup \{v\})$ contains an induced $P_3$, then 
    both $G-w$ and $(G\pivot vw)-w$ contain induced subgraphs isomorphic to $P_3+tP_1$, 
    contradicting the assumption that $G$ is $(P_3+tP_1)$-pivot-unique.
    So, $G-(A\cup \{v\})$ is $P_3$-free and it is a disjoint union of complete graphs.

    We claim that 
    \begin{itemize}
        \item[($\ast$)] $G$ has no two adjacent vertices that are twins in $G$.
    \end{itemize}
    Suppose that there are such vertices $x$ and $y$.
    By the assumption, $G$ contains a pivot-minor~$H$ isomorphic to $P_3+tP_1$.
    Since $P_3+tP_1$ has no two adjacent vertices that are twins, by Lemma~\ref{lem:truetwin}, $\{x,y\}\setminus V(H)$ is non-empty. Then by Lemma~\ref{lem:pivot-twin}, $G$ is not $(P_3+tP_1)$-pivot-unique, a contradiction. Therefore, the claim holds.

    We claim that $G[N_G(v)]- A$ contains no edges. Assume that there exists an edge $uw$ in $G[N_G(v)]- A$. Note that $u$ and $w$ are twins because every connected component of $G-(A\cup \{v\})$ is a complete graph. This contradicts $(\ast)$.
By a similar observation, $G[V(G)\setminus A\setminus N_G(v)]$ contains no edges.

    By the above claim, for every connected component $C$ of $G-(A\cup \{v\})$,
    each of $V(C)\setminus N_G(v)$ and $V(C)\cap N_G(v)$ contains at most one vertex.  
    Thus, either $C$ is a single vertex, or $C$ has two vertices where $v$ has only one neighbour in $C$.
    This shows that 
      $G$ is isomorphic to a graph obtained from $K_{1,n} + mP_1$ for some integers $n\ge t+1$ and $m\ge 0$, by 1-subdividing at most $n-(t+1)$ edges.
 \end{proof}

\begin{lemma}\label{lem:bounding}
  Let $G$ be a graph and $H$ be a pivot-minor of $G$.  
  Let $X$ be a set of vertices of degree at most $1$ in $G$
  and let $Y$ be the set of all neighbours of vertices in $X$. 
  If  $G-x$ does not have a pivot-minor isomorphic to $H$ for all $x\in X\cup Y$, 
  then $\abs{X\cup Y}\le \abs{V(H)}$.
\end{lemma}
\begin{proof}
  By the assumption on $X\cup Y$, $H$ contains every isolated vertex of~$G$.
  If $x\in X$ 
  and $y$ is its neighbour, 
  then $G\pivot xy$ is identical to the graph obtained from $G$ by swapping $x$ and $y$. 
  Thus none of $G-x$, $G\pivot xy-x$, $G-y$, and $G\pivot xy-y$ contains a pivot-minor isomorphic to~$H$. 
  Thus, $x,y\in V(H)$.
  Therefore $X\cup Y\subseteq V(H)$.
\end{proof}

\begin{lemma}\label{lem:P3iso2}
    Let $t$ be a positive integer and let $G$ be a graph obtained from $K_{1,n}+mP_1$ for some integers $n\ge t+1$ and $m\ge 0$ by $1$-subdividing at most $n-(t+1)$ edges. 
    \begin{enumerate}[label=\rm(\roman*)]
    \item If $G$ is an induced-subgraph-minimal graph having a pivot-minor isomorphic to $P_3+tP_1$, 
    then $\abs{V(G)}\le t+3$.
    \item Let $w$ be a vertex of degree $2$ in $G$, and let $v$ be the neighbour of $w$ having degree more than $1$ (corresponding to the center of $K_{1,n}$). If $G\pivot vw$ is an induced-subgraph-minimal graph having a pivot-minor isomorphic to $P_3+tP_1$, 
    then $\abs{V(G)}\le t+5$.
    \end{enumerate}
\end{lemma}

\begin{proof}
    Note that every vertex of degree at least $2$ is adjacent to a vertex of degree $1$.
    By \ref{lem:bounding}, $\abs{V(G)}\le t+3$, proving (i).

\medskip
       Now, we prove (ii). 
       Note that $G-v$ has no pivot-minor isomorphic to $P_3+tP_1$ because every connected component of $G-v$ has at most two vertices.
       Thus, by applying Lemma~\ref{lem:bounding} with the set $X$ of all vertices of degree at most $1$ non-adjacent to $v$,
       we deduce that $\abs{V(G)}-2\le t+3$.
  \end{proof}

\begin{lemma}\label{lem:pivotequivalent}
    Let $G$ be a graph obtained from $K_{1,n}+mP_1$ for some integers $n\ge 1$ and $m\ge 0$, by $1$-subdividing at least one edge. Let $vw$ be an edge of $G$ whose both ends have degree at least $2$. Then every graph pivot-equivalent to $G$ is isomorphic to $G$ or $G\pivot vw$.
\end{lemma}
\begin{proof}
    As pivoting any edge incident with a vertex of degree $1$ produces an isomorphic graph, every graph obtained from $G$ by one edge pivot is isomorphic to $G$ or $G\pivot vw$. It suffices to show that any graph obtained from $G\pivot vw$ by one edge pivot is again isomorphic to $G$ or $G\pivot vw$.

    Note that $v$ or $w$ has degree $2$ in $G$. Without loss of generality, we assume that $w$ has degree $2$, and 
    let $z$ be the neighbour of $w$ other than $v$ in $G$. Let $A$ be the set of vertices of degree $1$ in $G$ that are adjacent to $v$, and let $B\subseteq V(G)\setminus \{w\}$ be the set of vertices of degree $2$ in $G$ that are adjacent to $v$.
    
    The graph $G\pivot vw$ is obtained from $G$ by adding all edges between $z$ and $A\cup B$, and then swapping $v$ and $w$. Observe that the graphs obtained from $G\pivot vw$ by pivoting an edge between $A\cup \{v\}$ and $\{w,z\}$ are isomorphic to each other, and the graphs obtained from $G\pivot vw$ by pivoting an edge between $B$ and $\{w,z\}$ are isomorphic to each other. 
    It is not hard to verify that every graph obtained from $G\pivot vw$ by pivoting an edge between $A\cup \{v\}$ and $\{w,z\}$ is isomorphic to $G$, and every graph obtained from $G\pivot vw$ by pivoting an edge between $B$ and $\{w,z\}$ is isomorphic to $G\pivot vw$. 
    
    This proves the lemma.
\end{proof}

\begin{theorem}\label{thm:P3tP1}
    Let $t$ be a positive integer. If $G$ is an induced-subgraph-minimal graph having a pivot-minor isomorphic to $P_3+tP_1$, 
    then $\abs{V(G)}\le (t+1)(2^{t+3}-t-4)+2$.
\end{theorem}
\begin{proof}
    By Lemma~\ref{lem:pivotunique}, $G$ is $(P_3+tP_1)$-pivot-unique.
    Since $G$ has a pivot-minor isomorphic to $P_3+tP_1$, 
    there exists a graph~$H$ pivot-equivalent to $G$ such that $H$ contains an induced subgraph isomorphic to $P_3+tP_1$. 
    By Lemma~\ref{lem:pivot-pu}, $H$
    is $(P_3+tP_1)$-pivot-unique.
    Let $X\subseteq V(H)$ such that $H[X]$ is isomorphic to $P_3+tP_1$.
    For a subset $Y$ of $X$, let $N_Y$ be the set of vertices~$v$ in $V(H)\setminus X$ such that $N_H(v)\cap X=Y$.

    By Lemma~\ref{lem:pivot-twin}, 
    $N_Y$ is independent for all $Y\subseteq X$.
    By Lemma~\ref{lem:pivot-isolated}, $N_\emptyset=\emptyset$.

    Let $v_1v_2v_3$ be the induced path of $H[X]$.
    We claim that 
    \begin{itemize}
        \item[($\ast$)] for every $Y\subseteq X$, if $Y\setminus \{v_1,v_2,v_3\}\neq\emptyset$, then $\abs{N_Y}\le \min(\abs{Y}-1,1)$.    
    \end{itemize}
    Assume that $Y\subseteq X$ contains an isolated vertex of $H[X]$.
    If $\abs{Y}=1$, then by Lemma~\ref{lem:pivot-twin}, $N_Y=\emptyset$, proving the claim. 
    So we may assume that $\abs{Y}>1$ and $\abs{N_Y}>1$.
    Let $y$ be an isolated vertex of $H[X]$ in $Y$ and let $x\in N_Y$.
    Then in $H\pivot xy$, there are no edges between $X\setminus \{y\}$ and $N_Y\setminus \{x\}$, and so $H[(X\cup N_Y)\setminus \{x,y\}]$  contains an induced subgraph isomorphic to $P_3+tP_1$. Thus, both $H-x$ and $(H\pivot xy)-x$ contain pivot-minors isomorphic to $P_3+tP_1$, contradicting the fact that $H$ is $(P_3+tP_1)$-pivot-unique.
    Thus, the claim holds.

    By Lemma~\ref{lem:pivot-twin}, 
    $N_{\{v_1,v_2\}}=N_{\{v_2,v_3\}}=N_{\{v_1,v_2,v_3\}}=\emptyset$.
    We claim that 
    \begin{itemize}
        \item[($\ast\ast$)] for every $Y\subseteq \{v_1, v_3\}$, $\abs{N_Y}\le 1$.
    \end{itemize}
        Suppose for contradiction that there exists $Y\subseteq \{v_1, v_3\}$ such that $\abs{N_Y}\ge 2$.

    First assume that $Y=\{v_1,v_3\}$.
    Then both $H-v_1$ and $H\pivot v_1v_2-v_1$ contain induced subgraphs isomorphic to $P_3+tP_1$. Thus, we may assume that $\abs{Y}=1$ and $Y=\{v_i\}$ for some~$i\in \{1,3\}$. By symmetry, assume that $i=1$. Then $H-v_2$ and $H\pivot v_2v_3-v_2$ contain induced subgraphs isomorphic to $P_3+tP_1$.
    Therefore, the claim holds.

    By $(\ast)$ and $(\ast\ast)$, we have 
    \[
      \left\lvert\bigcup_{Y\subseteq X, Y\neq\{v_2\}} N_Y\right\rvert
    \le 
    (2^{t+3}-8-t)+3=2^{t+3}-t-5
    \]

    If $v_2$ has more than $t$ neighbours of degree~$1$ in $H$, 
    then by Lemma~\ref{lem:P3iso1}, $H$ is isomorphic to  a graph obtained from $K_{1,n} + mP_1$ for some integers $n\ge t+1$ and $m\ge 0$, by 1-subdividing at most $n-(t+1)$ edges. 
    Since $G$ is pivot-equivalent to $H$, by Lemma~\ref{lem:pivotequivalent}, $G$ is isomorphic to $H$ or $H\pivot vw$ where $vw$ is an edge whose both ends have degree at least $2$.
    In both cases, by Lemma~\ref{lem:P3iso2}, $G$ has at most $t+5$ vertices, and we are done. 
    Therefore, we may assume that $v_1$ has at most~$t$ neighbours of degree~$1$ in $H$. 

    Assume that 
    \(\abs{N_{\{v_2\}}}-t
    > t(2^{t+3}-t-5)\).
    Then 
    by the pigeonhole principle, 
    $H$ has a vertex~$q\in \bigcup_{Y\subseteq X,Y\neq \{v_2\}} N_Y $ such that $q$ has more than $t$ neighbors in $N_{\{v_2\}}$.
    Let $U=N_{\{v_2\}}\cap N_H(q)$
    and $z\in U$.
    Observe that in $H\pivot qz$, 
    $X$ is anti-complete to $U\setminus \{z\}$.   
    If $v_2$ is adjacent to any vertex $r$ in $X\setminus \{v_2\}$ in $H\pivot qz$, 
    then the induced path $qv_2r$ and the vertices in $U\setminus \{q\}$ form an induced $P_3+tP_1$. This means that both $H-z$ and $H\pivot qz-z $ contain pivot-minors isomorphic to~$P_3+tP_1$. So we deduce that $v_2$ has no neighbours in $X\setminus \{v_2\}$ in~$H\pivot qz$. 
    This implies that $q\in N_{\{v_1,v_3\}}$. Then $H\pivot qz-q$ contains a pivot-minor isomorphic to $P_3+tP_1$ on $(X\setminus \{v_2\})\cup \{z\}$. This contradicts the fact that $H$ is $(P_3+tP_1)$-pivot-unique.
  
    Therefore $\abs{N_{\{v_2\}}}
    \le t(2^{t+3}-t-5)+t$.
    Then $\abs{V(G)}\le (t+3)+(2^{t+3}-t-5)+t(2^{t+3}-t-5)+t\le (t+1)(2^{t+3}-t-4)+2$.
  \end{proof}

\begin{theorem}\label{t-general}
Let $t$ be a positive integer, and let $H=tP_1$, $H=P_2+tP_1$, $H=P_3+tP_1$, or $H=K_{1,t}$.
Then there is a polynomial-time algorithm for $H$-{\sc Pivot-Minor}  that gives an $H$-pivot-minor-sequence if one exists.
\end{theorem}

\begin{proof}
By Theorems~\ref{thm:tP1}, \ref{thm:P2tP1}, \ref{thm:star} and~\ref{thm:P3tP1}, the set ${\cal F}_H$ of induced-subgraph-minimal graphs containing a pivot-minor isomorphic to $H$ consists of finitely many non-isomorphic graphs.  
One can enumerate all non-isomorphic graphs in $\mathcal{F}_H$ in constant time, as $H$ is a fixed graph. By testing the existence of an induced subgraph isomorphic to a graph in $\mathcal{F}_H$ from a given graph $G$, one can in polynomial time either find an induced subgraph of $G$ that is isomorphic to a graph in $\mathcal{F}_H$ or correctly decides that $G$ does not contain a pivot-minor isomorphic to~$H$. When we find an induced subgraph $F$ of $G$ that is isomorphic to some graph in $\mathcal{F}_H$, in polynomial time, we can find the vertex deletions and edge pivots that modify $F$ into $H$. 
\end{proof}

\section{When $\mathbf{H}$ Is Fixed: $\mathbf{|V(H)|\le 4}$}\label{s-small}

We  give a certifying algorithm for recognizing $H$-pivot-minor-free graphs for every graph $H$ on at most four vertices except for the cases where $H \in \{K_4,C_3+P_1\}$.
Results in Section~\ref{s-tp1free} imply that there is such an algorithm for $H$ where
\[H\in \{P_1,2P_1,P_2,3P_1, P_1+\nobreak P_2, P_3, P_1+\nobreak P_2,2P_2,P_3,4P_1, P_1+\nobreak P_3, 2P_1+P_2, K_{1,3}\}.\]
Thus, it remains to show for $H$ where 
\[H\in \{C_3, 2P_2, P_4, C_4, \allowbreak \mbox{paw},\allowbreak \mbox{diamond}\}.\]For each such graph~$H$, we determine the set~${\cal F}_H$ of induced-subgraph-minimal graphs containing a pivot-minor isomorphic to~$H$.
The case $H=2P_2$ is the most involved, and we present the result at the end of this section. 

We first show that the class of $C_3$-pivot-minor-free graphs is exactly the class of bipartite graphs. 
\begin{proposition}\label{p-c3}
The following statements are equivalent for every graph~$G$:
\begin{enumerate}[label=\rm(\roman*)]
\item $G$ is $C_3$-pivot-minor-free.
\item $G$ has no induced odd cycle.
\item $G$ is bipartite.
\end{enumerate}
\end{proposition}

\begin{proof}
A graph is bipartite if and only if it has no odd cycle.
Note that for every odd integer~$n\ge 3$, $C_n$ contains a pivot-minor isomorphic to~$C_3$.
Hence, (i) implies (ii).
It is known that the class of bipartite graphs is pivot-minor-closed (see~\cite{Ou05}).
Hence~every bipartite graph does not contain a pivot-minor isomorphic to~$C_3$, and (iii) implies (i).
\end{proof}

Next, we consider $P_4$-pivot-minor-free graphs. As $P_4$ is pivot-equivalent to $C_4$, these graphs are also $C_4$-pivot-minor-free graphs. We show that these graphs are exactly the graphs whose connected components are obtained from stars by replacing each vertex with a clique. To describe this structure, we introduce clique-stars. 

A graph is a \emph{clique-star} if it is either a complete graph or it consists of pairwise vertex-disjoint cliques~$K$, $L_1,\ldots,L_p$ for some $p\geq 1$, such that every vertex of~$K$ is adjacent to every vertex of $L_1 \cup \cdots \cup L_p$ and there is no edge between any two distinct cliques~$L_i$ and~$L_j$. Note that every complete graph is also a clique-star with $p=1$.

\begin{lemma}\label{lem:clique-star}
The class of clique-stars is pivot-minor-closed.
\end{lemma}
\begin{proof}
Let $G$ be a clique-star. As every pivot-minor of a complete graph is a complete graph, we may assume that it consists of pairwise vertex-disjoint cliques $(K,L_1,\ldots,L_p)$ for some integer $p\geq 2$ such that every vertex of~$K$ is adjacent to every vertex of $L_1 \cup \cdots \cup L_p$ and there is no edge between any two distinct cliques~$L_i$ and~$L_j$. Let $uv\in E(G)$.

If~$u$ and~$v$ both belong to~$K$ or both belong to some~$L_i$, then pivoting~$uv$ results in the same graph.
In the remaining case, we may assume that~$u$ belongs to~$K$ and~$v$ belongs to~$L_1$.
After pivoting~$uv$, we obtain a new clique-star with $K'= L_1$, $L_1'= K$ and $L_i'=L_i$ for $i\in \{2,\ldots,p\}$.
Hence, the class of clique-stars is pivot-minor-closed.
\end{proof}

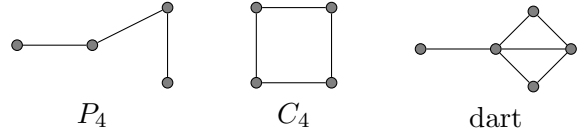
\begin{figure}
  \centering
  \begin{center}

    \begin{tikzpicture}
      \tikzstyle{v}=[circle,draw,fill=black!50,inner sep=0pt,minimum width=4pt]
      \draw (0,0) node[v](v){};
      \draw (1,.5) node[v](v1){};
      \draw(1,-.5)node[v](v2){};
      \draw (-1,0) node[v](v3){};
      \draw(v3)--(v)--(v1)--(v2);
      \draw(0,-.5) node [label=below:$P_4$]{};
    \end{tikzpicture}
    $\quad\quad$
    \begin{tikzpicture}
      \tikzstyle{v}=[circle,draw,fill=black!50,inner sep=0pt,minimum width=4pt]
      \draw (0,.5) node[v](v){};
      \draw (1,.5) node[v](v1){};
      \draw(1,-.5)node[v](v2){};
      \draw (0,-.5) node[v](v3){};
      \draw(v)--(v1)--(v2)--(v3)--(v);
      \draw(0.5,-.5) node [label=below:$C_4$]{};
    \end{tikzpicture}
    $\quad\quad$
    \begin{tikzpicture}
      \tikzstyle{v}=[circle,draw,fill=black!50,inner sep=0pt,minimum width=4pt]
      \draw (0,0) node[v](v){};
      \draw (.5,.5) node[v](v1){};
      \draw(.5,-.5)node[v](v2){};
      \draw (-1,0) node[v](v3){};
      \draw(1,0)node[v](w){};
      \draw(v3)--(v)--(v1)--(w)--(v2)--(v)--(w);
      \draw(0,-.5) node [label=below:dart]{};
    \end{tikzpicture}
  \end{center}
  \caption{The set~$\mathcal{F}_{P_4}$ of forbidden induced subgraphs for~$P_4$-pivot-minor-free graphs.}
\label{fig:nop4}
\end{figure}

We use Lemma~\ref{lem:clique-star} to prove the following result.

\begin{proposition}\label{p-p4}
The following statements are equivalent for every graph~$G$.
\begin{enumerate}[label=\rm(\roman*)]
\item $G$ is $P_4$-pivot-minor-free.
\item $G$ is $C_4$-pivot-minor-free.
\item $G$ is $(P_4, C_4, \text{dart})$-free
{\em (see \figurename~\ref{fig:nop4})}.
\item $G$ is the disjoint union of clique-stars.
\end{enumerate}
\end{proposition}

\begin{proof}
Both the $P_4$ and~$C_4$ can be obtained from each other by pivoting one edge and so (i) and (ii) are equivalent.

We observe that the dart contains~$P_4$ as a pivot-minor. Let $v$ be a vertex of degree $2$ and $w$ be the vertex of degree $3$ in the dart $H$. Then $(H\pivot vw)[V(H)\setminus \{w\}]$ is isomorphic to $P_4$. Thus, (i) implies (iii).

Lemma~\ref{lem:clique-star} implies that the class of graphs all of whose connected components are clique-stars is pivot-minor-closed, hence
(iv) implies (i). 

It remains to prove that (iii) implies (iv).
Suppose that~$G$ has a connected component~$D$ that is not a clique-star.
Also assume that $G$ is $(P_4,C_4)$-free.
It is well known that  the complement of a connected $P_4$-free graph on at least two vertices is disconnected~\cite{CLB81}.
Hence, we can partition~$V(D)$ into two non-empty sets~$A$ and~$B$, such that $A$ is complete to $B$.
Moreover, as~$D$ is not a complete graph, we may assume that~$B$ is not a clique.
If~$A$ is not a clique either, then two non-adjacent vertices of~$A$, together with two non-adjacent vertices of~$B$, form an induced~$C_4$, a contradiction.
Hence~$A$ is a clique.
We may assume that~$A$ is chosen to be maximal subject to the condition that $A$ is complete to~$B$ and~$B$ has two non-adjacent vertices.

Suppose~$G[B]$ is connected. Again since $G$ is $P_4$-free and $G[B]$ has at least two vertices, we can partition~$B$ into two non-empty sets~$B_1$ and~$B_2$, such that $B_1$ is complete to~$B_2$.

As~$B$ is not a clique, this means that at least one of~$B_1$ and~$B_2$, say~$B_2$, is not a clique.
Then, by the same argument as before, $B_1$ must be a clique.
This implies that $A\cup B_1$ is complete to~$B_2$.
This contradicts the maximality of~$A$, as we could have chosen $A\cup B_1$ instead.
Hence~$G[B]$ is not connected.

Let $J_1,\ldots,J_r$ be the connected components of~$G[B]$ for some $r\geq 2$.
If all of $J_1, \ldots, J_r$ are cliques, then $G[D]$ is a clique-star, contradicting the assumption that $G[D]$ is not a clique-star.
Thus, one of $J_1,\ldots,J_r$, say~$J_1$, is not a clique. Then $J_1$ contains an induced path $uvw$.
Then $u,v,w$, together with a vertex of~$A$ and a vertex of~$J_2$, induce a dart.

This shows that (iii) implies (iv).
\end{proof}

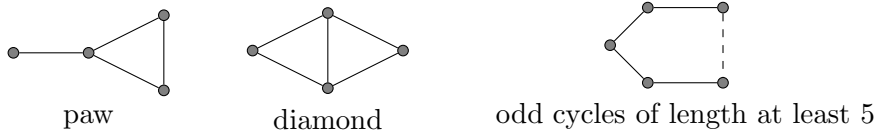
\begin{figure}
  \centering
  \begin{center}
    \begin{tikzpicture}
      \tikzstyle{v}=[circle,draw,fill=black!50,inner sep=0pt,minimum width=4pt]
      \draw (0,0) node[v](v){};
      \draw (1,.5) node[v](v1){};
      \draw(1,-.5)node[v](v2){};
      \draw (-1,0) node[v](v3){};
      \draw(v3)--(v)--(v1)--(v2)--(v);
      \draw(0,-.5) node [label=below:paw]{};
    \end{tikzpicture}
    $\quad\quad$
    \begin{tikzpicture}
      \tikzstyle{v}=[circle,draw,fill=black!50,inner sep=0pt,minimum width=4pt]
      \draw (0,0) node[v](w2){};
      \draw (1,.5) node[v](v1){};
      \draw(1,-.5)node[v](v2){};
      \draw (2,0) node[v](w1){};
      \foreach \i in {1,2}{
        \draw (w1)--(v\i)--(w2);
      }
      \draw (v1)--(v2);
      \draw(1,-.5) node [label=below:diamond]{};
    \end{tikzpicture}
    $\quad\quad$
    \begin{tikzpicture}
      \tikzstyle{v}=[circle,draw,fill=black!50,inner sep=0pt,minimum width=4pt]
      \draw (0,.5) node[v](v){};
      \draw (1,.5) node[v](v1){};
      \draw(1,-.5)node[v](v2){};
      \draw (0,-.5) node[v](v3){};
      \draw (-.5,0) node[v](w){};
      \draw (w)--(v)--(v1);
      \draw (v2)--(v3)--(w);
      \draw [dashed] (v1)--(v2);
      \draw(0.5,-.5) node [label=below:odd cycles of length at least~$5$]{};
    \end{tikzpicture}
  \end{center}
  \caption{The set $\mathcal{F}_{\text{paw}}$ of forbidden induced subgraphs for $\text{paw}$-pivot-minor-free graphs.}
\label{fig:nopaw}
\end{figure}

We now consider paw-pivot-minor-free graphs. As paw is pivot-equivalent to diamond, these graphs are also diamond-pivot-minor-free graphs. We show that these graphs are exactly graphs whose connected components are either bipartite or complete.

\begin{proposition}\label{p-paw}
The following statements are equivalent for every graph~$G$.
\begin{enumerate}[label=\rm(\roman*)]
\item $G$ is paw-pivot-minor-free.
\item $G$ is diamond-pivot-minor-free.
\item $G$ has no induced subgraph isomorphic to paw, the diamond or an odd cycle of length at least~$5$. 
{\em (see \figurename~\ref{fig:nopaw})}.
\item Every connected component of $G$ is either bipartite or complete.
\end{enumerate}
\end{proposition}

\begin{proof}
By pivoting one edge, the diamond can be obtained from the paw and so (i) and (ii) are equivalent.
Since every odd cycle on at least five vertices contains a pivot-minor isomorphic to the paw, (i) implies (iii).
As the class of graphs whose connected components are complete graphs or bipartite graphs is pivot-minor-closed, (iv) implies (i).

It remains to prove that (iii) implies (iv).
Suppose (iii) holds.
Let~$D$ be a connected component of~$G$.

We claim that~$D$ is bipartite or complete.
Assume that it is not true. As $D$ is not bipartite, it contains an induced odd cycle. As $D$ has no induced odd cycle of length at least $5$, $D$ has a triangle.
Let~$K$ be a maximal clique of~$D$. By the existence of a triangle, $K$ has size at least $3$.

As~$D$ is not complete and~$K$ is maximal, there is a vertex $u\in V(D)\setminus K$ that has both a neighbour and a non-neighbour in~$K$.
If $u$ has exactly one neighbour in~$K$, then it has at least two non-neighbours in~$K$, and this implies that $D$ contains an induced subgraph isomorphic to the paw.
If $u$ has at least two neighbours in~$K$, then $D$ contains an induced subgraph isomorphic to the diamond. These lead to a contradiction.
\end{proof}

Lastly, we characterize $2P_2$-pivot-minor-free graphs.  It turns out that there are nine induced-subgraph-minimal graphs having a pivot-minor isomorphic to $2P_2$ (see Figure~\ref{fig:2P_2}). Note that a $2P_2$-pivot-minor-free graph cannot have two connected components containing an edge. We show that $2P_2$-pivot-minor-free graphs containing an edge are exactly the graphs, in which the unique connected component containing an edge is either an induced subgraph of 
the prism or $W_5$ (see also Figure~\ref{fig:s2})
or a leaf-attached complete multipartite graph, which we define below.

\begin{figure}
\begin{center}
\begin{tikzpicture}
  \tikzstyle{v}=[circle,draw,fill=black!50,inner sep=0pt,minimum width=4pt]
  \draw (0,0) node [v] (v1){};
  \draw (1,0) node [v] (v2){};
  \draw (1,1) node [v] (v3){};
  \draw (0,1) node [v] (v4){};
    \draw (v2)--(v3);
    \draw (v4)--(v1);
     \draw (0.5,0) node[label=below:{$O_1=2P_2$}]{};
\end{tikzpicture}
$\quad\quad$
         \begin{tikzpicture}[rotate=0]
          \tikzstyle{v}=[circle,draw,fill=black!50,inner sep=0pt,minimum width=4pt]
      \draw (0,0) node[v](w2){};
      \draw (1,.5) node[v](v1){};
      \draw(1,-.5)node[v](v2){};
      \draw (2,0) node[v](w1){};
      \draw (3,0) node[v](z){};
      \foreach \i in {1,2}{
        \draw (w1)--(v\i)--(w2);
      }
      \draw (z)--(w1);
    \draw (1.5,-0.5) node[label=below:$O_2$]{};
    \end{tikzpicture}
$\quad\quad$
         \begin{tikzpicture}[rotate=0]
          \tikzstyle{v}=[circle,draw,fill=black!50,inner sep=0pt,minimum width=4pt]
      \draw (0,0) node[v](w2){};
      \draw (1,.5) node[v](v1){};
      \draw(1,-.5)node[v](v2){};
      \draw (2,0) node[v](w1){};
      \draw (3,0) node[v](z){};
      \foreach \i in {1,2}{
        \draw (w1)--(v\i)--(w2);
      }
      \draw (v1)--(v2);
      \draw (z)--(w1);
    \draw (1.5,-.5) node[label=below:$O_3$]{};
    \end{tikzpicture}
$\quad\quad$
\begin{tikzpicture}
  \tikzstyle{v}=[circle,draw,fill=black!50,inner sep=0pt,minimum width=4pt]
  \draw (0,0) node [v] (v1){};
  \draw (0:1) node [v] (v2){};
  \draw (0:2) node [v] (v3){};
  \draw (60:1) node [v] (w2){};
  \draw (60:2) node [v] (w3){};
  \draw (30:1.73205080757) node [v] (x1){};
  \draw (v1)--(v2)--(v3)--(x1)--(w3)--(w2)--(v1);
  \draw (v2)--(w2)--(x1)--(v2);
    \draw (1,0) node[label=below:$O_4$]{};
\end{tikzpicture}

\bigskip
\noindent
\begin{tikzpicture}
  \tikzstyle{v}=[circle,draw,fill=black!50,inner sep=0pt,minimum width=4pt]
  \draw (90:1) node [v] (v0){};
  \draw (0,0) node [v] (v){};
  \draw (270+45:1) node [v] (v1){};
  \draw (270+15:1) node [v] (v2){};
  \draw (270-15:1) node [v] (v3){};
  \draw (270-45:1) node [v] (v4) {};
  \foreach \i in {0,1,2,3,4} {
    \draw (v)--(v\i);
    }
    \draw (v1)--(v2)--(v3)--(v4);
    \draw (0,-1) node[label=below:$O_5$]{};
\end{tikzpicture}
$\quad\quad$
\begin{tikzpicture}[rotate=90]
  \tikzstyle{v}=[circle,draw,fill=black!50,inner sep=0pt,minimum width=4pt]
  \draw (0,0) node [v] (v00){};
  \draw (1,0) node [v] (v10){};
  \draw (2,0) node [v] (v20){};
  \draw (0,1) node [v] (v01){};
  \draw (1,1) node [v] (v11){};
  \draw (2,1) node [v] (v21){};
  \draw (v10) -- (v20) -- (v11) -- (v21) -- (v10) -- (v00) -- (v01) -- (v11) -- (v10);
    \draw (0,0.5) node[label=below:$O_6$]{};
\end{tikzpicture}
$\quad\quad$
\begin{tikzpicture}
  \tikzstyle{v}=[circle,draw,fill=black!50,inner sep=0pt,minimum width=4pt]
  \draw (0,0) node [v] (v1){};
  \draw (1,0) node [v] (v2){};
  \draw (1,1) node [v] (v3){};
  \draw (0,1) node [v] (v4){};
  \draw (.5,1.86602540378) node[v] (v) {};
  \foreach \i in {3,4} {
    \draw (v)--(v\i);
    }
    \draw (v1)--(v2)--(v3)--(v4)--(v1)--(v3);
    \draw (v2)--(v4);
    \draw (0.5,0) node[label=below:$O_7$]{};
\end{tikzpicture}
$\quad\quad$
\begin{tikzpicture}[rotate=90]
  \tikzstyle{v}=[circle,draw,fill=black!50,inner sep=0pt,minimum width=4pt]
  \draw (0,0) node [v] (v00){};
  \draw (1,0) node [v] (v10){};
  \draw (2,0) node [v] (v20){};
  \draw (0,1) node [v] (v01){};
  \draw (1,1) node [v] (v11){};
  \draw (2,1) node [v] (v21){};
  \draw (v10) -- (v20) -- (v11) -- (v21) -- (v10) -- (v00) -- (v01) -- (v11) -- (v10) -- (v01);
    \draw (0,0.5) node[label=below:$O_8$]{};
\end{tikzpicture}
$\quad\quad$
   \begin{tikzpicture}
    \tikzstyle{v}=[circle,draw,fill=black!50,inner sep=0pt,minimum width=4pt]
      \draw (0,0) node[v](w2){};
      \draw (1,.5) node[v](v1){};
      \draw(1,-.5)node[v](v2){};
      \draw (2,0) node[v](w1){};
      \draw (1,1.5) node[v](z){};
      \foreach \i in {1,2}{
        \draw (w1)--(v\i)--(w2);
      }
      \draw (v1)--(v2);
      \draw (w2)--(z)--(w1);
     \draw (1,-.5) node[label=below:$O_9$]{};
   \end{tikzpicture}
\end{center}
\caption{The set~${\cal F}_{2P_2}$ of nine minimal forbidden induced subgraphs for $2P_2$-pivot-minor-free graphs.}
\label{fig:2P_2}
\end{figure}
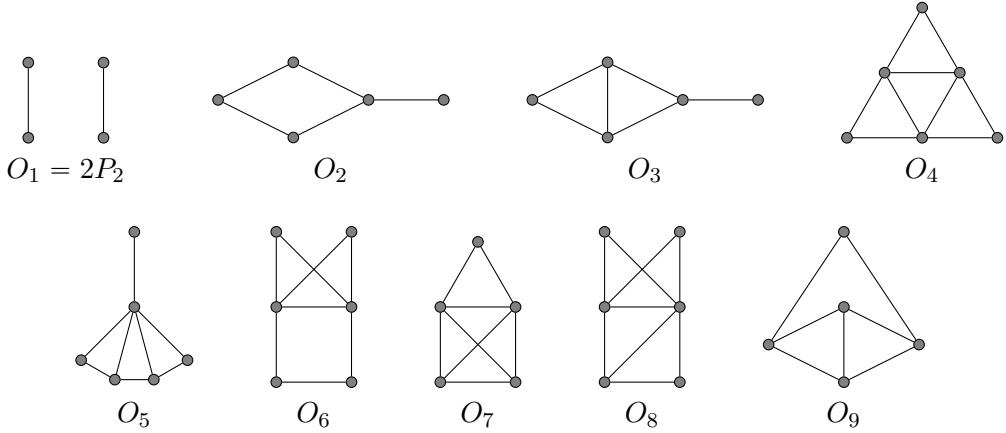

\begin{figure}[t]
 \centering
  \begin{tikzpicture}
    \tikzstyle{v}=[circle,draw,fill=black!50,inner sep=0pt,minimum width=4pt]
    \draw (0,2) node [v] (a1){};
    \draw (0.5,2) node [v] (a2){};
    \draw (1,2) node [v] (a3){};

    \draw (-1,1) node [v] (b1){};

    \draw (-2,1.25) node [v] (x1){};
    \draw (-2,1) node [v] (x2){};
    \draw (-2,0.75) node [v] (x3){};

    \draw (2,1.25) node [v] (c1){};
    \draw (2,0.75) node [v] (c2){};

    \draw (-0.5,0) node [v] (d1){};
    \draw (0,-0.25) node [v] (d2){};
    
    \draw (1.5,0) node [v] (e1){};

    \draw (1.5,-1) node [v] (y1){};
    \draw (1.75,-1) node [v] (y2){};

\foreach \i in {1,2,3}{
    \foreach \j in {1}{
    \draw (a\i)--(b\j);}}

\foreach \i in {1,2,3}{
    \foreach \j in {1,2}{
    \draw (a\i)--(c\j);}}

\foreach \i in {1,2,3}{
    \foreach \j in {1,2}{
    \draw (a\i)--(d\j);}}

\foreach \i in {1,2,3}{
    \foreach \j in {1}{
    \draw (a\i)--(e\j);}}

\foreach \i in {1}{
    \foreach \j in {1,2}{
    \draw (b\i)--(c\j);}}

\foreach \i in {1}{
    \foreach \j in {1,2}{
    \draw (b\i)--(d\j);}}

\foreach \i in {1}{
    \foreach \j in {1}{
    \draw (b\i)--(e\j);}}

\foreach \i in {1,2}{
    \foreach \j in {1,2}{
    \draw (c\i)--(d\j);}}

\foreach \i in {1,2}{
    \foreach \j in {1}{
    \draw (c\i)--(e\j);}}

\foreach \i in {1,2}{
    \foreach \j in {1}{
    \draw (d\i)--(e\j);}}

\foreach \i in {1,2,3}{
    \foreach \j in {1}{
    \draw (x\i)--(b\j);}}

\foreach \i in {1,2}{
    \foreach \j in {1}{
    \draw (y\i)--(e\j);}}
  \end{tikzpicture}
\caption{An example of a leaf-attached complete multipartite graph.}\label{fig:leafattached}
\end{figure}
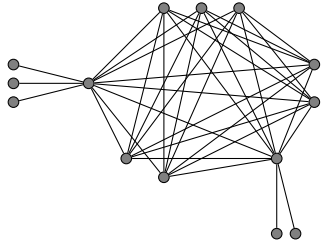

A graph $G=(V,E)$ is {\it complete multipartite} if  $V$ can be partitioned into $p\geq 1$ independent sets $V_1,\ldots, V_p$ called the {\it partition classes} of $G$, such that for any two vertices $u$ and $v$ it holds that $uv\in E$ if and only if $u\in V_i$ and $v\in V_j$ for some distinct $i, j\in \{1, \ldots, p\}$. If some $V_i$ has size~$1$, we say that the (unique) vertex of $V_i$ is a {\it singleton vertex} of $G$. We say that we {\it attach a leaf} to a vertex $v$ in a graph $G$ if we add a new vertex $u$ to $G$ and make $u$ only adjacent to $v$.
A \emph{leaf-attached complete multipartite graph} is a graph obtained from a complete multipartite graph~$G$ by attaching zero or more
leaves to each singleton vertex of $G$; 
see also Figure~\ref{fig:leafattached}.

\begin{lemma}\label{lem:multipartite}
    Every leaf-attached complete multipartite graph is $2P_2$-pivot-minor-free.
\end{lemma}

\begin{proof}
    Let $G$ be a leaf-attached complete multipartite graph. Say, $G$ is obtained from a complete multipartite graph~$G'$ by attaching leaves to singleton vertices of $G'$. 
    We prove the statement of the lemma by induction on $\abs{V(G)}$. 
    If $G'$ has only one partition class, then either $G$ is edgeless or $G$ is a star. So, $G$ has no pivot-minor isomorphic to $2P_2$. Now assume that $G'$ has at least two partition classes. This implies that $G'$ is connected, and thus $G$ is connected as well. It can be readily checked that if $\abs{V(G)}=4$, then the statement is true. 

    Assume that  $\abs{V(G)}\ge 5$ and suppose for contradiction that $G$ contains a pivot-minor isomorphic to $2P_2$.
     Then by Lemma~\ref{lem:bouchet}, there exists a vertex $w$ in $G$ such that $G-w$ or $G/w$ contains a pivot-minor isomorphic to $2P_2$.   
 If $w$ is a leaf, then the connected component of  $G-w$ or $G/w$ containing an edge is again a leaf-attached complete multipartite graph. 
 This would contradict our induction hypothesis.
 Thus, we may assume that $w$ is not a leaf. Furthermore, we may assume that there is no leaf adjacent to $w$; otherwise, we can apply the same argument. 

    Let $z$ be a neighbour of $w$. From the above, we conclude that $z$ is contained in another partition class of $G'$. Observe that $G-w$ and $(G\pivot wz)-w$ are the disjoint unions of a leaf-attached complete multipartite graph and isolated vertices. This again contradicts the induction hypothesis.
\end{proof}

\begin{lemma}\label{lem:hom}
  Let $G$ be an $O_2$-free graph.
  If two vertices $a$ and $b$ are adjacent and $N_{G}(a)\setminus(N_G(b)\cup \{b\}), N_G(b)\setminus(N_G(a)\cup \{a\})$ are independent sets, then 
  $N_{G}(a)\setminus(N_G(b)\cup \{b\})$ is anti-complete or complete to $N_G(b)\setminus(N_G(a)\cup \{a\})$.
\end{lemma}
\begin{proof}
  Suppose that this is not true. Then there exist vertices $x,y,z$ where 
  \begin{itemize}
      \item $\{x,z\}$ and $\{y\}$ are contained in  distinct sets of $N_{G}(a)\setminus(N_G(b)\cup \{b\})$ and $N_G(b)\setminus(N_G(a)\cup \{a\})$, 
      \item $xy\in E(G)$ and $yz\notin E(G)$.
  \end{itemize}
  Then $G[\{a,b,x,y,z\}]$ is isomorphic to $O_2$, a contradiction.  %
\end{proof}
\begin{proposition}\label{p-2P2}
The following statements are equivalent for every graph~$G$.
\begin{enumerate}[label=\rm(\roman*)]
\item $G$ is $2P_2$-pivot-minor free.
\item $G$ is $(O_1, \ldots, O_9)$-free (see Figure~\ref{fig:2P_2}).
\item There is at most one connected component of $G$ containing an edge, and if such a connected component $H$ exists, then either 
    \begin{itemize}
        \item $H$ is isomorphic to an induced subgraph of the prism or $W_5$, or 
        \item $H$ is a leaf-attached complete multipartite graph.
    \end{itemize}
\end{enumerate}
\end{proposition}
\begin{proof}
    It is readily seen that (i) implies (ii). 
    Moreover, (iii) implies (i) by Lemma~\ref{lem:multipartite} and the fact that the prism and $W_5$ are $2P_2$-pivot-minor-free.
  Hence, it remains to show that (ii) implies (iii), which we do below.
  
  Suppose that $G$ is an $(O_1, \ldots, O_9)$-free graph.
    We may assume that $G$ contains an edge, and that $G$ has no isolated vertices. As $G$ is $2P_2$-free, this means that $G$ is connected.

  Let $C$ be a maximum clique in $G$. We may assume $|C|\ge 2$. 

First assume that $|C|=2$. 
    First consider the case when for every edge $vw$, one of $N_G(w)\setminus \{v\}$ and $N_G(v)\setminus \{w\}$ is empty. If there is an edge $vw$ where both sets are empty, then $G$ is isomorphic to $P_2$. Otherwise, $G$ is isomorphic to a star, so $G$ is a complete multipartite graph. Thus, we may assume that there is an edge $vw$, where both $N_G(w)\setminus \{v\}$ and $N_G(v)\setminus \{w\}$ are non-empty.
  As the size of a maximum clique in $G$ is $2$, we can let $C=\{v,w\}$. 
      
    For each $Y\subseteq C$, let $A_Y$ be the set of vertices $v$ in $V(G)\setminus C$ with $N_G(v)\cap C=Y$. As $C$ is a maximum clique, $A_C=\emptyset$, and $A_{\{v\}}$ and $A_{\{w\}}$ are independent sets. 
By our choice of~$v$ and~$w$, we have that both $A_{\{v\}}$ and $A_{\{w\}}$ are non-empty.
 As $G$ is $2P_2$-free, we also find that $A_\emptyset$ is an independent set. 

 By Lemma~\ref{lem:hom}, $A_{\{v\}}$ is anti-complete or complete to $A_{\{w\}}$.  
If $A_\emptyset=\emptyset$, then this implies that $G$ is a leaf-attached complete multipartite graph, irrespective of whether $A_{\{v\}}$ and $A_{\{w\}}$ are anti-complete or complete to each other, as $A_{\{v\}}$ and $A_{\{w\}}$ are independent sets. Hence, we may assume that $A_\emptyset\neq \emptyset$.

    As $G$ is connected 
     and $A_\emptyset$ is an independent set, every vertex in $A_\emptyset$ has a neighbour in $A_{\{v\}}\cup A_{\{w\}}$. 
First suppose that $A_{\{v\}}$ is complete to $A_{\{w\}}$. Since $G$ is $O_2$-free, $A_\emptyset$ is complete to $A_{\{v\}}\cup A_{\{w\}}$. However, as $A_\emptyset$ is non-empty, this contradicts our assumption that $C$ is a maximum clique.
Hence, $A_{\{v\}}$ is anti-complete to $A_{\{w\}}$. As $G$ is $2P_2$-free, $A_\emptyset$ must be complete to $A_{\{v\}}\cup A_{\{w\}}$. As $G$ is $O_2$-free, each of $A_{\{v\}}$ and $A_{\{w\}}$ has size $1$. As $G$ is $O_2$-free and $A_\emptyset$ is non-empty, $A_\emptyset$ has size~$1$. Consequently, $G$ is isomorphic to $C_5$, which is an induced subgraph of $W_5$.

Now assume that $|C|\ge 3$.
    Let $C=\{v_1, \ldots, v_m\}$ for some integer~$m\geq 3$.
    We first show that every $w\in N_G(C)$ has either exactly one neighbour in $C$ or exactly one non-neighbour in $C$. Since $C$ is a maximum clique, $w$ has at least one non-neighbour in $C$. If $w$ has at least two neighbours in $C$ and at least two non-neighbours in $C$, then $G$ contains an induced subgraph isomorphic to $O_7$. Thus, the claim holds.

    For every $i\in \{1, \ldots, m\}$, let $A_i$ be the set of vertices $v$ in $V(G)\setminus C$ such that $N_G(v)\cap C=\{v_i\}$, and let $B_i$ be the set of vertices $v$ in $V(G)\setminus C$ such that $N_G(v)\cap C=C\setminus \{v_i\}$. We observe that each $A_i$ is independent as $G$ is $2P_2$-free and $m\geq 3$, and that each $B_i$ is independent as $C$ is a maximum clique. 

    We claim that every vertex in $V(G)\setminus C$ has a neighbour in $C$. For a contradiction, suppose that there is a vertex $v$ in $V(G)\setminus C$ having no neighbour in $C$. As $G$ is connected and it is $2P_2$-free, $v$ has a neighbour $w\in N_G(C)$. If $N_G(w)\cap C=\{v_i\}$ for some $i\in\{1,2,\ldots,m\}$, then $vw$ and an edge in $C\setminus \{v_i\}$ form an induced $2P_2$, and if $N_G(w)\cap C=C\setminus \{v_i\}$ for some $i\in\{1,2,\ldots,m\}$, then $C\cup \{v,w\}$ contains an induced $O_3$. In both cases we obtain a contradiction. Thus, every vertex in $V(G)\setminus C$ has a neighbour in $C$.

    Next, we claim that either $A_i=\emptyset$ or $B_i=\emptyset$ for every $i\in \{1, \ldots, m\}$. Assume that for some $i\in \{1, \ldots, m\}$, $A_i$ and $B_i$ are non-empty. Let $v\in A_i$ and $w\in B_i$. If $vw\in E(G)$, then $C\cup \{v,w\}$ contains an induced $O_9$, and otherwise, $C\cup \{v,w\}$ contains an induced $O_3$. In both cases we obtain a contradiction. Thus, the claim holds.

    Let $I$ be the set of indices $i$ such that $A_i$ is non-empty, and let $J$ be the set of indices $j$ such that $B_j$ is non-empty.  By the above observation, we find that $I\cap J=\emptyset$.

   We now distinguish between the following three cases.

    \medskip
\noindent
\textbf{(Case 1. $G$ contains an edge between  $A_{i_1}$ and $A_{i_2}$ for some $i_1,i_2\in I$.)}
\nopagebreak

Recall that every $A_i$ is an independent set. Hence, $i_1\neq i_2$.    
We may assume without loss of generality that $\{1,2\}\subseteq I$, and that $x\in A_1$ and $y\in A_2$ are adjacent. If $\abs{C}\ge 4$, then $G[\{v_3, v_4, x, y\}]$ is isomorphic to $2P_2$, a contradiction. Thus, $\abs{C}=3$.
As $I\cap J=\emptyset$ and $A_1$ and $A_2$ are non-empty, we find that $B_1=B_2=\emptyset$. Hence, $V(G)=\{v_1,v_2,v_3\}\cup A_1\cup A_2\cup A_3\cup B_3$, where either $A_3=\emptyset$ or $B_3=\emptyset$.

    First, we show that $A_1=\{x\}$ and $A_2=\{y\}$. 
    If $y$ has another neighbour $x'$ in $A_1$, then $G[\{v_1, x, x', y, v_3\}]$ is isomorphic to $O_2$, a contradiction. By symmetry, $x$ has no other neighbour in $A_2$. As $A_1$ and $A_2$ are independent sets, this means that the connected component of $G[A_1\cup A_2]$ that contains $x$ and $y$ only consists of $x$ and $y$. 
    Since $G$ is $2P_2$-free, $G[A_1\cup A_2]$ has no two connected components containing an edge.
    If there is an isolated vertex $x'$ in $G[A_1\cup A_2]$, then $G[\{v_1, v_2, x, y, x'\}]$ is isomorphic to $O_2$, a contradiction. This shows that $A_1=\{x\}$ and $A_2=\{y\}$.  

We now show that $B_3=\emptyset$. Suppose $B_3$ contains a vertex $z$. If $z$ is anti-complete to $\{x,y\}$, then $G[\{v_1, v_2, v_3, x,y,z\}]$ is isomorphic to $O_6$. If $z$ is complete to $\{x,y\}$, then $G[\{v_2, v_3, x,y,z\}]$ is isomorphic to $O_3$. If $z$ has exactly one neighbour in $\{x,y\}$, then $G[\{v_1, v_2, x,y,z\}]$ is isomorphic to $O_9$. Each of these cases is a contradiction, and we conclude that $B_3=\emptyset$.

    Lastly, we prove that every vertex of $A_3$ is complete to $\{x,y\}$. Let $z\in A_3$. If $z$ is anti-complete to $\{x,y\}$, then $G[\{x,y, z, v_3\}]$ is isomorphic to $2P_2$, a contradiction. If $z$ has exactly one neighbour in $\{x,y\}$, then $G[\{v_1, v_2, x,y,z\}]$ is isomorphic to $O_2$, another contradiction. We conclude that $z$ is complete to $\{x,y\}$.
Now applying the first claim to the pair $(A_1, A_3)$ instead of $(A_1,A_2)$ yields that $A_3$ has size at most $1$. Thus, $G$ is isomorphic to an induced subgraph of the prism.

    \medskip
\noindent
\textbf{(Case 2. $G$ contains a non-edge  between $B_{j_1}$ and $B_{j_2}$ for some $j_1,j_2\in J$ with $j_1\neq j_2$.)}
          
   We may assume without loss of generality that $\{1,2\}\subseteq J$, and that $x\in B_1$ and $y\in B_2$ are non-adjacent.
If $\abs{C}\ge 4$, then $G[\{v_1, v_3, v_4, x, y\}]$ is isomorphic to $O_7$, a contradiction. Thus, $\abs{C}=3$.
 As $I\cap J=\emptyset$ and $B_1$ and $B_2$ are non-empty, we find that $A_1=A_2=\emptyset$. Hence, $V(G)=\{v_1,v_2,v_3\}\cup A_3\cup B_1\cup B_2\cup B_3$, where either $A_3=\emptyset$ or $B_3=\emptyset$.

    First, we show that $B_1=\{x\}$ and $B_2=\{y\}$. 
     Recall that $B_1$ and $B_2$ are independent sets.   
   Hence, if $y$ has another non-neighbour $x'$ in $B_1$, then $G[\{v_1, v_2, v_3, x, x', y\}]$ is isomorphic to~$O_8$,    
    a contradiction. By symmetry, we find that $x$ has no other non-neighbours in $B_2$. If $y$ has some neighbour $x'$ in $B_1$, then $G[\{v_1, v_2, x, x', y\}]$ is isomorphic to $O_2$, another contradiction. By symmetry, we conclude that $B_1=\{x\}$ and $B_2=\{y\}$.

     Second, we show that $\abs{A_3}\le 1$ and that $A_3$ is complete to $\{x,y\}$. To start with showing the latter, let $z\in A_3$. If $z$ is anti-complete to $\{x,y\}$, then $G[\{v_1, v_2, v_3, x,y,z\}]$ is isomorphic to $O_5$, a contradiction. If $z$ is adjacent to exactly one of $x,y$, say $x$, then $G[\{x,y,v_1,z\}]$ is isomorphic to $2P_2$, another contradiction. Thus, $z$ is complete to $\{x,y\}$. Now, if $A_3$ has two vertices $z, z'$, then $G[\{x,y,z,z',v_1\}]$ is isomorphic to $O_2$, again a contradiction. Thus, $\abs{A_3}\le 1$.

We first  assume that $\abs{A_3}= 1$, say $A_3=\{z\}$ for some vertex~$z$, so $B_3=\emptyset$. As $A_3$ is complete to $\{x,y\}$, we now find that $G=G[\{v_1,v_2,v_3,x,y,z\}]$ is isomorphic to $W_5$, where $v_3$ is the vertex of degree $5$. In the remainder, we assume that $A_3=\emptyset$. We show that $B_3=\emptyset$. For a contradiction, suppose that there exists a vertex $q\in B_3$. If $q$ is anti-complete to $\{x,y\}$, then $G[\{v_1, v_2, v_3, x, y, q\}]$ is isomorphic to $O_4$. If $q$ has exactly one neighbour in $\{x,y\}$, say $y$, then $G[\{v_1, v_2, x, y, q\}]$ is isomorphic to $O_3$. If $q$ is complete to $\{x,y\}$, then 
     $G[\{v_2, v_3, x, y, q\}]$ is isomorphic to $O_9$. So, $B_3=\emptyset$.
     This implies that $G$ is isomorphic to the gem, which is an induced subgraph of $W_5$.

   \medskip
\noindent
\textbf{(Case 3. For all distinct $i_1, i_2\in I$, $A_{i_1}$ is anti-complete to $A_{i_2}$, and for all distinct $j_1, j_2\in J$, $B_{j_1}$ is complete to $B_{j_2}$.)}

    First assume that $\abs{C}\ge 4$. 
    We claim that $\bigcup_{i\in I}A_i$ is anti-complete to $\bigcup_{j\in J}B_j$. To see this, assume that there are $i\in I$, $j\in J$, and distinct $r, r'\in \{1, \ldots, m\}\setminus \{i,j\}$ with $x\in A_i$ and $y\in B_j$ such that $x$ is adjacent to $y$. Then $G[\{x, y, v_j, v_{r}, v_{r'}\}]$ is isomorphic to $O_3$, a contradiction. Thus, the claim holds.
   Therefore, $G$ is
     a leaf-attached complete multipartite graph.

From now on we assume that $\abs{C}=3$.  
We claim that 
   for all $i\in I$ and $j\in J$, 
   $A_i$ is anti-complete or complete to $B_j$. 
   To prove this, 
   suppose that $A_i\cup B_j$ has three vertices $x$, $x'$, $y$ such that 
   $xy\in E(G)$, $x'y\notin E(G)$, and either $A_i\cap \{x,x',y\}=\{y\}$ or $B_j\cap \{x,x',y\}=\{y\}$.
   We recall that every $A_i$ and every $B_j$ is an independent set and so $x$ is non-adjacent to $x'$.
    Now, let $i\in I$ and $j\in J$ and recall that $I\cap J=\emptyset$, so $i\neq j$. 
    Then $y\notin B_j$ because otherwise $x,x'\in A_i$ and 
    $G[\{y,x,x'\}\cup C]$ is isomorphic to $O_5$.
    We deduce that $y\in A_i$ and $x,x'\in B_j$ and therefore 
    $G[\{x,y,y'\}\cup C]$ is isomorphic to $O_8$, contradicting our assumption.
    Therefore $A_i$ is anti-complete or complete to~$B_j$.
     
We claim that 
    if $A_i$ is complete to $B_j$ for some $i\in I$ and $j\in J$, then $\abs{A_i}=1$ and $\abs{B_j}=1$. If $\abs{A_i}\ge 2$, then $G[A_i\cup B_j\cup C]$ has an induced subgraph isomorphic to $O_8$, a contradiction.
    If $\abs{B_j}\ge 2$, then $G[A_i\cup B_j\cup (C\setminus \{v_i\})]$ has an induced subgraph  isomorphic to $O_2$, another contradiction. So, indeed we have that $\abs{A_i}=1$ and $\abs{B_j}=1$.
    
   If $I=\emptyset$ or $J=\emptyset$, then
    $G$ is a leaf-attached complete multipartite graph.
    Thus, we may assume that $I\neq \emptyset$ and $J\neq\emptyset$. As $I\cap J=\emptyset$ and $\abs{C}=3$, we have that $2\le  \abs{I}+\abs{J}\le 3$. We discuss each of these three cases separately. Recall that for every $i\in I$ and $j\in J$, it holds that $A_i$ is complete or anti-complete to $B_j$.

   First assume that $\abs{I}=1$ and $\abs{J}=1$, say $I=\{1\}$ and $J=\{2\}$. If $A_1$ is anti-complete to $B_2$, then $G$ is a leaf-attached complete multipartite graph. If $A_1$ is complete to $B_2$, then $\abs{A_1}=1$ and $\abs{B_2}=1$. This implies that $G$ is isomorphic to the gem, which is an induced subgraph of $W_5$.

   Now assume that $\abs{I}=1$ and $\abs{J}=2$, say $I=\{1\}$ and thus $J=\{2, 3\}$. 
If $A_1$ is anti-complete to $B_2\cup B_3$, then $G$ is a leaf-attached complete multipartite graph.
   Otherwise we may assume without loss of generality that $A_1$ is complete to $B_2$. So, $\abs{A_1}=1$ and $\abs{B_2}=1$. 
   If $A_1$ is anti-complete to $B_3$, then $G- v_1$ contains an induced subgraph isomorphic to $O_2$, a contradiction. 
    So, $A_1$ is also complete to $B_3$. 
   Then $\abs{B_3}=1$.
   However, now $A_1\cup B_2\cup B_3 \cup \{v_1\}$ is a clique on four vertices, contradicting the maximality of $C$.

   Finally, assume that $\abs{I}=2$ and $\abs{J}=1$, say $I=\{1, 2\}$ and $J=\{3\}$. 
If $B_3$ is anti-complete to $A_1\cup A_2$, then $G$ is a leaf-attached complete multipartite graph.   
   Otherwise
   we may assume without loss of generality that $A_1$ is complete to $B_3$. So, $\abs{A_1}=1$ and $\abs{B_3}=1$. If $A_2$ is complete to $B_3$, then $\abs{A_2}=1$ and $G$ is isomorphic to $O_4$, a contradiction. Thus, $A_2$ is anti-complete to $B_3$. 
   We now find that $G[A_1\cup A_2\cup B_3\cup\{v_1,v_2\}]$ contains an induced subgraph isomorphic to $O_3$, a contradiction.
This completes the proof of the theorem.
\end{proof}

\begin{theorem}\label{t-1}
For $H\in \{C_3, 2P_2, P_4, C_4, \allowbreak \mbox{paw},\allowbreak \mbox{diamond}\}$,
there is a polynomial-time algorithm for $H$-{\sc Pivot-Minor} that gives an $H$-pivot-minor-sequence if one exists.
\end{theorem}

\begin{proof}
If $H=C_3$, by Proposition~\ref{p-c3}, we need to find an odd cycle~$F$, which we do in polynomial time by testing bipartiteness.

If $H\in\{P_4, C_4\}$, then we use condition~(iii) in Proposition~\ref{p-p4} to decide if a graph has a pivot-minor isomorphic to~$H$ and obtain an induced subgraph $F$ of $G$ that is isomorphic to a graph in~${\cal F}_H$ if one exists.

Assume $H\in\{\mbox{paw},\mbox{diamond}\}$. 
Let $G$ be a given graph. For each connected component~$D$ of~$G$, we test the bipartiteness of $D$. If it is bipartite, then we skip to other connected component. Otherwise, we find in polynomial time an odd cycle $C$ in $D$. If $C$ has length at least $5$, then by condition (iii) in Proposition~\ref{p-paw}, it is a graph in ${\cal F}_H$. Assume that $C$ has length $3$. Then we find a maximal clique $C^*$ containing $C$. If there is a vertex $V(D)\setminus C^*$, then there is a vertex in $V(D)\setminus C^*$ having a neighbour and a non-neighbour in $C^*.$ This provides an induced subgraph~$F$ isomorphic to either the paw or the diamond, which are graphs in ${\cal F}_H$.

If $H=2P_2$, then we use condition~(ii) in Proposition~\ref{p-2P2} to decide if a graph has a pivot-minor isomorphic to~$H$ and obtain an induced subgraph $F$ of $G$ that is isomorphic to a graph in~${\cal F}_H$ if one exists.

The theorem follows, as in polynomial time we can find the vertex deletions and edge pivots
that modify $F$ into $H$.
\end{proof}

\section{The Two Open Cases}\label{s-infi}

\begin{figure}[b]
  \centering
  \begin{tikzpicture}
    \tikzstyle{v}=[circle,draw,fill=black!50,inner sep=0pt,minimum width=4pt]
      \draw (-1,.5) node[v](lv){};
      \draw (-1,-.5) node[v](lv3){};
      \draw (-.5,0) node[v](w){};
      \draw (0,.5) node[v](v){};
      \draw (0,-.5) node[v](v3){};
      \draw (w)--(lv)--(lv3)--(w)--(v)--(v3)--(w);
      \node at (-0.5,-.5)[label=below:$\overline{C_4+P_1}$]{};
  \end{tikzpicture}
  $\quad\quad$
  \begin{tikzpicture}
    \tikzstyle{v}=[circle,draw,fill=black!50,inner sep=0pt,minimum width=4pt]
    \begin{scope}[xshift=-2cm]
      \draw (0,.5) node[v](lv){};
      \draw (-1,.5) node[v](lv1){};
      \draw(-1,-.5)node[v](lv2){};
      \draw (0,-.5) node[v](lv3){};
      \draw (.5,0) node[v](lw){};
      \draw (-.5,0) node {odd};
\end{scope}
      \draw (lw)--(lv)--(lv1);
      \draw (lv2)--(lv3)--(lw);
      \draw [dashed] (lv1)--(lv2);
      \draw (0,.5) node[v](v){};
      \draw (1,.5) node[v](v1){};
      \draw(1,-.5)node[v](v2){};
      \draw (0,-.5) node[v](v3){};
      \draw (-.5,0) node[v](w){};
      \draw (0,0) node {odd};
      \draw (w)--(v)--(v1);
      \draw (v2)--(v3)--(w);
      \draw [dashed] (v1)--(v2);
      \draw [dashed](w)--(lw);
      \node at (0,-.5)[label=below:Two odd induced cycles joined by a path]{};
  \end{tikzpicture}
  \caption{An infinite family of minimal forbidden induced subgraphs for the class of $K_4$-pivot-minor-free graphs.}
  \label{fig:nok4}
\end{figure}
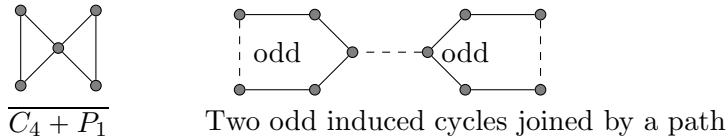

We recall that it is still an open problem to solve {\sc $H$-Pivot-Minor} if $H\in \{K_4,C_3+P_1\}$. In this brief section, we only show that the sets~${\cal F}_{K_4}$ and~${\cal F}_{C_3+P_1}$ consist of infinitely many non-isomorphic graphs.

\begin{proposition}
The set~${\cal F}_{K_4}$ has infinitely many non-isomorphic graphs.
\end{proposition}

\begin{proof}
The set~${\cal F}_{K_4}$ contains the following infinite class of graphs.
Let~${\cal H}$ be the class of graphs obtained from taking two odd cycles and adding a path of length~$0$ or more connecting one vertex of one
odd cycle with one vertex of the other
odd cycle (see \figurename~\ref{fig:nok4}).
Note that no graph of~${\cal H}$ is isomorphic to an induced subgraph of another graph of~${\cal H}$.
Moreover, every graph in~${\cal H}$ contains a pivot-minor isomorphic to the graph~$\overline{C_4+P_1}$, while~$\overline{C_4+P_1}$ contains a pivot-minor isomorphic to~$K_4$.
Finally, the graphs in~${\cal H}$ can be shown to be minimal, as any graph obtained from identifying a vertex of some odd cycle
with an end-vertex of a path that ends in the centre of a claw does not contain~$K_4$ as a pivot-minor.
\end{proof}

\begin{proposition}
The set~${\cal F}_{C_3+P_1}$ has infinitely many non-isomorphic graphs.
\end{proposition}

\begin{proof}
As mentioned in the proof of Proposition~\ref{p-c3}, the class of bipartite graphs is pivot-minor-closed (see also~\cite{Ou05}).
Hence, no graph in~${\cal F}_{C_3+P_1}$ is bipartite.
For every odd integer $k\geq 3$, the graph $C_k+\nobreak P_1$ contains~$C_3+\nobreak P_1$ as a pivot-minor.
Therefore~${\cal F}_{C_3+P_1}$ contains all graphs of this form.
\end{proof}

\section{Determining the Set $\mathbf{{\cal F}_{K_{1,3}}}$}\label{s-claw}

\begin{figure}
  \begin{center}
  \begin{tikzpicture}
    \tikzstyle{v}=[circle,draw,fill=black!50,inner sep=0pt,minimum width=4pt]
    \draw (0,0) node [v] (v1){};
    \draw (1,0) node [v] (v2){};
    \draw (1,1) node [v] (v3){};
    \draw (0,1) node [v] (v){};
    \foreach \i in {1,2,3}{
      \draw (v)--(v\i);
    }
    \draw (0.5,0) node[label=below:$K_{1,3}$]{};
  \end{tikzpicture}
  $\quad\quad$
  \begin{tikzpicture}
    \tikzstyle{v}=[circle,draw,fill=black!50,inner sep=0pt,minimum width=4pt]
    \draw (0,0) node [v] (v1){};
    \draw (1,0) node [v] (v2){};
    \draw (1,1) node [v] (v3){};
    \draw (0,1) node [v] (v4){};
    \draw (.5,.5) node[v] (v) {};
    \draw (v4)--(v1)--(v)--(v2)--(v3);
    \draw (0.5,0) node[label=below:$P_5$]{};
  \end{tikzpicture}
  $\quad\quad$
  \begin{tikzpicture}
    \tikzstyle{v}=[circle,draw,fill=black!50,inner sep=0pt,minimum width=4pt]
    \draw (0,0) node [v] (v1){};
    \draw (1,0) node [v] (v2){};
    \draw (1,1) node [v] (v3){};
    \draw (0,1) node [v] (v4){};
    \draw (.5,.5) node[v] (v) {};
    \draw (v4)--(v1)--(v)--(v2)--(v3);
    \draw (v1)--(v2);
    \draw (0.5,0) node[label=below:bull]{};
  \end{tikzpicture}
  $\quad\quad$
  \begin{tikzpicture}
    \tikzstyle{v}=[circle,draw,fill=black!50,inner sep=0pt,minimum width=4pt]
    \draw (0,0) node [v] (v1){};
    \draw (1,0) node [v] (v2){};
    \draw (1,1) node [v] (v3){};
    \draw (0,1) node [v] (v4){};
    \draw (.5,.5) node[v] (v) {};
    \foreach \i in {1,2,3,4} {
      \draw (v)--(v\i);
      }
      \draw (v1)--(v2)--(v3)--(v4)--(v1);
      \draw (0.5,0) node[label=below:$W_4$]{};
  \end{tikzpicture}
  $\quad\quad$
    \begin{tikzpicture}
      \tikzstyle{v}=[circle,draw,fill=black!50,inner sep=0pt,minimum width=4pt]
    \begin{scope}[xshift=-1.5cm]
      \foreach \i in {1,2} {
      \draw (120*\i:.7) node[v,label=left:$a_\i$](a\i){};
    }
      \foreach \i in {3} {
      \draw (-10:.7) node[v,label=left:$a_\i$](a\i){};
    }
    \draw (a1)--(a2)--(a3)--(a1);
    \end{scope}
    \begin{scope}[xshift=1.5cm]
      \draw (10:1) node[v,label=right:$b_4$] (b4){};
      \foreach \i in {1,2} {
      \draw (180-120*\i:.7) node[v,label=right:$b_\i$](b\i){};
      \draw (a\i)--(b\i)--(b4);
    }
      \foreach \i in {3} {
      \draw (190:.7) node[v,label=above left:$b_\i$](b\i){};
      \draw (a\i)--(b\i)--(b4);
    }
    \draw (b1)--(b2)--(b3)--(b1);
    \end{scope}
    \draw (-90:.4) node[label=below:$\bw_3$]{};
  \end{tikzpicture}
  \end{center}
\caption{The five graphs of the set ${\cal F}_{K_{1,3}}$.}
\label{fig:3P1etc-graphs}
\end{figure}
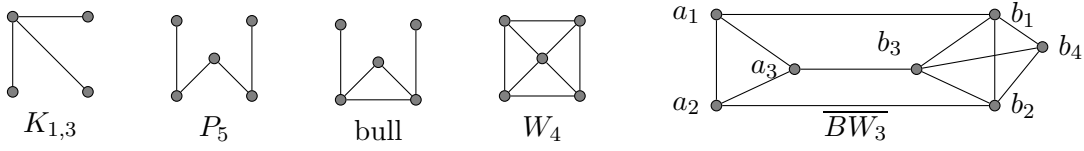

In this section, we first determine $\mathcal{F}_{3P_1}$ exactly and then use this result to determine $\mathcal{F}_{K_{1,3}}$. The latter set is displayed in \figurename~\ref{fig:3P1etc-graphs}.

\begin{theorem}\label{t-claw}
The following two statements hold:
\begin{enumerate}[label=\rm(\roman*)]
\item A graph $G$ is $3P_1$-pivot-minor-free if and only if it is $(3P_1, W_4, \bw_3)$-free.
\item A graph $G$ is $K_{1,3}$-pivot-minor-free if and only if it is $(K_{1,3}, P_5, \text{bull}, W_4, \bw_3)$-free.
\end{enumerate}
\end{theorem}
We prove in Lemma~\ref{lem:3P_1-bullclawP5-correspondence}
that a graph~$G$ is $(\text{\textnormal{bull}},K_{1,3},P_5)$-free if and only if every connected component of~$G$ is $3P_1$-free. Using this, the statement (ii) will follow from (i). Thus, we first show the statement (i). 
That is, in Proposition~\ref{p-3P_1}, we prove by induction on $|V(G)|$ that if a graph~$G$ contains $3P_1$ as a pivot-minor isomorphic, then $G$ contains a graph from $\{3P_1, W_4, \bw_3\}$ as an induced subgraph (the reverse implication is immediate).

The above claim holds if $|V(G)| \leq 3$, and so we may assume that $|V(G)|\ge 4$.
If~$G$ has a pivot-minor isomorphic to~$3P_1$, then by Lemma~\ref{lem:bouchet}, there is a vertex $v\in V(G)$ such that~$G-v$ or~$G/v$ contains a pivot-minor isomorphic to $3P_1$.
By the induction hypothesis, we may assume that~$G/v$ contains a pivot-minor isomorphic to~$3P_1$.
In Lemmas~\ref{lem:subcase1}, \ref{lem:subcase2} and~\ref{lem:subcase3}, we will show that if~$G/v$ contains an induced subgraph isomorphic to~$3P_1$, $W_4$ or~$\bw_3$, then so does~$G$. These lemmas will form the main steps in our induction.
We start by proving Lemma~\ref{lem:easy2}, which deals with some special cases.

\begin{figure}
  \begin{center}
  \begin{tikzpicture}
    \tikzstyle{v}=[circle,draw,fill=black!50,inner sep=0pt,minimum width=4pt]
    \tikzstyle{b}=[circle,draw,inner sep=0pt,minimum width=30pt]
      \draw (-1,1) node [v,label=$v$] (v){};
    \draw (1,1) node [v,label=$w$] (w){};
    \draw (0,0.2) node [b] (c){$S_2$};
    \draw (-2,-.7) node[b] (l){$S_1$};
    \draw (2,-.7) node[b] (r){$S_3$};
    \draw (0,-1.4)node[b] (x) {$S_4$};
    \draw (v)--(w)--(c)--(v);
    \draw (v)--(l);
    \draw (w)--(r);
    \draw [snake=zigzag](l)--(r)--(c)--(l);
    \draw[snake=zigzag] (l)--(x)--(r);
    \draw[snake=zigzag] (c)--(x);
  \end{tikzpicture}
  \end{center}
\caption{The sets $S_1, S_2, S_3, S_4$ in~$G$.}
\label{fig:pivot}
\end{figure}
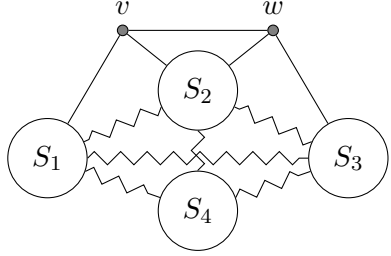

\begin{lemma}\label{lem:easy2}
Let~$vw$ be an edge of a graph~$G$, and let Let $S_1:=N_G(v)\setminus (N_G(w)\cup \{w\})$, $S_2:=N_G(v)\cap N_G(w)$, $S_3:=N_G(w)\setminus (N_G(v)\cup \{v\})$ and $S_4:=V(G)\setminus (N_G(v)\cup N_G(w)\cup \{v,w\})$
 (see \figurename~\ref{fig:pivot}).
If~$G\pivot vw$ satisfies at least one of the following conditions, then~$G$ contains an induced subgraph isomorphic to~$3P_1$ or~$W_4$:
\begin{enumerate}[label=\rm(\roman*)]
\item\label{lem:easy2-i}
There exist $v_1\in S_1$, $v_2\in S_2$ and $v_3\in S_3$ such that $\{v_1,v_2,v_3\}$ is a clique or an independent set in~$G\pivot vw$.
\item\label{lem:easy2-ii}
There exist distinct $v_1, v_2\in S_2$ and $v_3\in S_1\cup S_3$ such that~$v_1$ and~$v_2$ are non-adjacent and~$v_3$ is complete or anti-complete to $\{v_1, v_2\}$ in~$G\pivot vw$.
\item\label{lem:easy2-iii}
There exists an induced path $v_1v_2v_3$ in~$G\pivot vw$ and a vertex $v_4\in S_4$ such that~$v_4$ is complete to~$\{v_1, v_2, v_3\}$ in~$G\pivot vw$, and $\{v_1,v_2,v_3\}\subseteq S_1\cup S_2$ or $\{v_1,v_2,v_3\}\subseteq S_2\cup S_3$.
\item\label{lem:easy2-0} $S_1\cup S_4$ or $S_3\cup S_4$ is not a clique in~$G\pivot vw$.
\item\label{lem:easy2-v} $G\pivot vw$ has an induced cycle~$C$ of length~$4$ such that $V(C)\subseteq S_1\cup S_2\cup \{w\}$ or $V(C)\subseteq S_2\cup S_3\cup \{w\}$.
\end{enumerate}
\end{lemma}

\begin{proof}
Statements (i), (ii), (iii) are trivial.

For (iv), assume that $S_1\cup S_4$ is not a clique in~$G\pivot vw$.
Note that $(G\pivot vw)[S_1\cup S_4]=G[S_1\cup S_4]$.
Two non-adjacent vertices in $S_1\cup S_4$ and~$w$ form an independent set of size~$3$ in~$G$. We can prove in the same way when $S_3\cup S_4$ is not a clique in $G\pivot vw$.

Lastly, suppose that $G\pivot vw$ satisfies (v). We first assume that $V(C)\subseteq S_1\cup S_2\cup \{w\}$. The proof for the case $V(C)\subseteq S_2\cup S_3\cup \{w\}$ would be symmetric.

By (iv), we may assume that $S_1$ is a clique. This implies that $\abs{V(C)\cap S_1}\le 2$.  

First assume that $\abs{V(C)\cap S_1}=0$. Then $V(C)\subseteq S_2\cup \{w\}$. Since $w$ is complete to $S_2$ in $G\pivot vw$, $C$ does not contain $w$ and thus $V(C)\subseteq S_2$. This implies that Then $V(C)\cup \{v\}$ induces $W_4$ in $G$.

Assume that $\abs{V(C)\cap S_1}=1$. Then the two neighbours in $C$ of the vertex in $V(C)\cap S_1$ are contained in $S_2$, and the result follows by (ii).

Lastly, we assume that $\abs{V(C)\cap S_1}= 2$.
Let $C=abcda$ and without loss of generality, we assume that $a, b\in S_1$. As $c$ and $d$ are not complete to $\{a,b\}$ in $G\pivot vw$, we have $\{c, d\}\subseteq S_2$. 
In $G$, $acdba$ is an induced cycle dominated by $v$.
Therefore $V(C)\cup\{v\}$ induces~$W_4$ in~$G$.
\end{proof}

\begin{lemma}\label{lem:subcase1}
Let~$vw$ be an edge of a graph~$G$.
If~$(G\pivot vw)- v$ contains an induced subgraph isomorphic to~$3P_1$, then~$G$ contains an induced subgraph isomorphic to~$3P_1$ or~$W_4$.
\end{lemma}

\begin{proof}
Let $S_1:=N_G(v)\setminus (N_G(w)\cup \{w\})$, $S_2:=N_G(v)\cap N_G(w)$, $S_3:=N_G(w)\setminus (N_G(v)\cup \{v\})$ and $S_4:=V(G)\setminus (N_G(v)\cup N_G(w)\cup \{v,w\})$.

Let~$T$ be an independent set of size~$3$ in~$(G\pivot vw)-v$.
By Lemma~\ref{lem:easy2}\ref{lem:easy2-0}, we may assume that $S_1\cup S_4$ and $S_3 \cup S_4$ are cliques in $G\pivot vw$.
If $w\in T$, then the other two vertices are contained in $S_3\cup S_4$ which is a clique. This is a contradiction.
Thus, $w\notin T$. 

If~$\abs{T\cap S_2}=1$, then $\abs{T\cap S_1}=\abs{T\cap S_3}=1$ and by Lemma~\ref{lem:easy2}\ref{lem:easy2-i}, we are done in this case.
If $T\subseteq S_2 \cup S_4$, then~$T$ is independent in~$G$.
So we may assume that $|T \cap (S_1 \cup S_3)|=1$ and $|T\cap S_2|=2$. The proof is completed by Lemma~\ref{lem:easy2}\ref{lem:easy2-ii}.
\end{proof}

\begin{lemma}\label{lem:subcase2}
Let~$vw$ be an edge of a graph~$G$.
If $(G\pivot vw)-v$ contains an induced subgraph isomorphic to~$W_4$, then~$G$ contains an induced subgraph isomorphic to~$3P_1$, $W_4$ or~$\bw_3$.
\end{lemma}

\begin{proof}
Let $S_1:=N_G(v)\setminus (N_G(w)\cup \{w\})$, $S_2:=N_G(v)\cap N_G(w)$, $S_3:=N_G(w)\setminus (N_G(v)\cup \{v\})$ and $S_4:=V(G)\setminus (N_G(v)\cup N_G(w)\cup \{v,w\})$.

We proceed by induction on~$\abs{V(G)}$.
Let $H$ be a pivot-minor of $(G\pivot vw)-v$ that is isomorphic to $W_4$.
If~$V(G)\setminus V(H)\setminus \{v,w\}$ contains a vertex $x$ that is not in the~$W_4$, then we are done by induction, because $(G\pivot vw)-v-x=((G-x)\pivot vw)-v$.
Therefore we may assume that $(G\pivot vw)-v=H$ or $(G\pivot vw)-v-w=H$.
Let~$z$ be the vertex of degree $4$ in~$H$ and let $C=v_1v_2v_3v_4v_1$ be the induced cycle of length $4$ in~$H$.

Suppose that~$G$ has no induced subgraph isomorphic to~$3P_1$ or~$W_4$. By Lemma~\ref{lem:easy2}\ref{lem:easy2-0}, both $S_1\cup S_4$ and $S_3\cup S_4$ are cliques in~$G\pivot vw$.

\medskip
\noindent
\textbf{(Case 1. $(G\pivot vw)-v= H$.)}

If $z=w$, then $V(C)\subseteq S_1\cup S_2$. Then by Lemma~\ref{lem:easy2}\ref{lem:easy2-v}, $G$ contains an induced $3P_1$ or $W_4$. So, we may assume  that $z\neq w$.

Without loss of generality, we may assume that $v_1=w$.
Then $v_4zv_2$ is an induced path in $S_1\cup S_2$ in $G\pivot vw$ and so $v_3\notin S_4$ by Lemma~\ref{lem:easy2}\ref{lem:easy2-iii}.
As $v_1$ and $v_3$ are not adjacent in $G\pivot vw$, we have $v_3\in S_3$.
By Lemma~\ref{lem:easy2}\ref{lem:easy2-ii} and \ref{lem:easy2-v}, $v_4\in S_1$ or $v_2\in S_1$.
As~$S_1$ is a clique, we have either
\begin{itemize}
    \item $v_4\in S_1$ and $v_2\in S_2$, or
    \item $v_2\in S_1$ and $v_4\in S_2$.
\end{itemize} This contradicts Lemma~\ref{lem:easy2}\ref{lem:easy2-i}, because $z\in S_1\cup S_2$.

\medskip
\noindent
\textbf{(Case 2. $(G\pivot vw)-v-w= H$.)}

We may assume that~$G\pivot vw$ is non-isomorphic to~$G$ and therefore $S_1\cup S_3\neq\emptyset$.

Assume that $v_1,v_2,v_3,v_4\in S_2\cup S_4$. Since $S_1\cup S_3\neq \emptyset$, we have $z\in S_1\cup S_3$.
By Lemma~\ref{lem:easy2-ii}, $S_2$ is a clique.
As $S_4$ is also a clique, without loss of generality, we may assume that $v_1,v_2\in S_2$ and $v_3,v_4\in S_4$.
In this case~$G$ is isomorphic to~$\bw_3$.
By the symmetry of~$S_1$ and~$S_3$ we may therefore assume that $\{v_1,v_2,v_3,v_4\}\cap S_1\neq\emptyset$.

Assume $z\in S_4$. Then by Lemma~\ref{lem:easy2}\ref{lem:easy2-iii}, $\abs{S_1}+\abs{S_2}\le 2$ and $\abs{S_2}+\abs{S_3}\le 2$.
Then $\abs{S_2}\le 1$ because $S_1\neq\emptyset$.
If $\abs{S_1\cup S_2\cup S_3}=4$, then $\abs{S_1}+\abs{S_2}+\abs{S_3}=4$ and so $\abs{S_1}=\abs{S_3}=2$ and $S_2=\emptyset$.
As~$S_1$ and~$S_3$ are cliques, without loss of generality, we may assume that $v_1,v_2\in S_1$ and $v_3, v_4\in S_3$.
Then $G-v-w$ is isomorphic to~$W_4$, a contradiction.
If $\abs{S_1\cup S_2\cup S_3}=3$, then we may assume that $v_1\in S_4$.
Since $S_1\cup S_4$ and $S_3\cup S_4$ are cliques
$\abs{S_1}\le 1$ and $\abs{S_3}\le 1$.
Therefore $\abs{S_1}=\abs{S_2}=\abs{S_3}=1$.
Then~$G$ isomorphic to~$\bw_3$.
So we may assume that $z\notin S_4$.

Assume $S_4\neq\emptyset$. Without loss of generality,  we assume that $v_1\in S_4$.
Since $S_1\cup S_4$ and $S_3\cup S_4$ are cliques, $v_3\in S_2$. Recall that $\{v_1,v_2,v_3,v_4\}\cap S_1\neq\emptyset$.
By symmetry between~$v_2$ and~$v_4$, we may assume $v_2\in S_1$.
Then $z \notin S_3$ by Lemma~\ref{lem:easy2}\ref{lem:easy2-i}.
Since $S_1\cup S_4$ is a clique, $v_4\in S_2\cup S_3$.
By Lemma~\ref{lem:easy2}\ref{lem:easy2-iii}, $v_4\notin S_2$ and so $v_4\in S_3$.
By Lemma~\ref{lem:easy2}\ref{lem:easy2-i}, $z \notin S_1$, so $z \in S_2$.
Then~$G$ is isomorphic to~$\bw_3$.
So we may assume that $S_4=\emptyset$.
By Lemma~\ref{lem:easy2}\ref{lem:easy2-v}, $\{v_1,v_2,v_3,v_4\}\cap S_3\neq\emptyset$.

Lastly, assume $z\in S_2$. Then we may assume $v_1\in S_1$.
Then by Lemma~\ref{lem:easy2}\ref{lem:easy2-i}, $v_2,v_4\notin S_3$ and $v_3\in S_3$.
By Lemma~\ref{lem:easy2}\ref{lem:easy2-i}, $v_2,v_4\notin S_1$ and so $v_2,v_4\in S_2$, contradicting Lemma~\ref{lem:easy2}\ref{lem:easy2-ii}.
So we may assume that $z\in S_1\cup S_3$.
By symmetry between~$S_1$ and~$S_3$, we may assume that $z\in S_1$.
If $v_1\in S_2$, then by Lemma~\ref{lem:easy2}\ref{lem:easy2-i}, $v_2,v_4\notin S_3$ and so $v_3\in S_3$.
By Lemma~\ref{lem:easy2}\ref{lem:easy2-i} again, $v_2,v_4\notin S_2$ and therefore $v_2,v_4\in S_1$, contradicting the assumption that~$S_1$ is a clique.
So we may assume that $v_1,v_2,v_3,v_4\notin S_2$.
Since~$S_1$ and~$S_3$ are cliques, we may assume that $v_1,v_2\in S_1$ and $v_3,v_4\in S_3$.
Then~$G$ is isomorphic to~$\bw_3$.
This proves the lemma.
\end{proof}

\begin{lemma}\label{lem:subcase3}
Let~$G$ be a graph containing an edge~$vw$.
If~$(G\pivot vw)-v$ contains an induced subgraph isomorphic to~$\bw_3$, then~$G$ contains an induced subgraph isomorphic to~$3P_1$, $W_4$ or~$\bw_3$.
\end{lemma}
\begin{proof}
Let $S_1:=N_G(v)\setminus (N_G(w)\cup \{w\})$, $S_2:=N_G(v)\cap N_G(w)$, $S_3:=N_G(w)\setminus (N_G(v)\cup \{v\})$ and $S_4:=V(G)\setminus (N_G(v)\cup N_G(w)\cup \{v,w\})$.

Again, we proceed by the induction on~$\abs{V(G)}$.
Let $H$ be a pivot-minor of $(G\pivot vw)-v$ that is isomorphic to $\bw_3$.
If~$V(G)\setminus V(H)\setminus \{v,w\}$ contains a vertex $x$ that is not in the~$W_4$, then we are done by induction, because $(G\pivot vw)-v-x=((G-x)\pivot vw)-v$.
Thus, we may assume that $(G\pivot vw)-v$ or $(G\pivot vw)-v-w$ is isomorphic to $\bw_3$.
Let $U_1=\{a_1,a_2,a_3\}$ and $U_2=\{b_1,b_2,b_3,b_4\}$ be two disjoint cliques in~$H$ such that $a_1b_1,a_2b_2,a_3b_3$ are the three edges between $U_1$ and $U_2$ in~$H$.

Suppose~$G$ has no induced subgraph isomorphic to~$3P_1$, $W_4$ or~$\bw_3$.
By Lemma~\ref{lem:easy2}\ref{lem:easy2-0}, $S_1\cup S_4$ and $S_3\cup S_4$ are cliques of~$G\pivot vw$. We may assume that $S_1\cup S_3\neq\emptyset$ because otherwise $G\pivot vw=G$.

We prove the following claims (A)--(H).
\begin{enumerate}[label=\rm(\Alph*)]
\item Let $J\subseteq \{1,2,3\}$ of size $2$. If $\{a_j:j\in J\}\subseteq S_2$, then the vertex of $\{a_j:j\in \{1,2,3\}\setminus J\}$ is not contained in $S_1\cup S_3$.
\item Let $J\subseteq \{1,2,3\}$ of size $2$. If $\{b_j:j\in J\}\subseteq S_2$, then the vertices of $\{b_j:j\in \{1,2,3,4\}\setminus J\}$ are not contained in $S_1\cup S_3$.
\end{enumerate}
Suppose (A) does not hold.
By symmetry, we assume that $J=\{1,2\}$. We may assume that $a_3\in S_1$.
Then $b_1,b_2\notin S_1\cup S_4$ because $S_1\cup S_4$ is a clique.
Then $a_1a_2b_2b_1a_1$ is an induced cycle of $G\pivot vw$ contained in $\{w\}\cup S_2\cup S_3$, contradicting Lemma~\ref{lem:easy2}\ref{lem:easy2-v}. Thus, (A) holds, and similarly, (B) also holds.

\begin{enumerate}[label=\rm(\Alph*)]
\setcounter{enumi}{2}
\item Let $J\subseteq \{1,2,3\}$ of size $2$. If $\{a_j:j\in J\}\subseteq S_4$ and $\{b_j:j\in J\}\subseteq S_i$ for some $i=1,2,3$, then the vertex of $\{a_j:j\in \{1,2,3\}\setminus J\}$ is contained in $S_i\cup S_4\cup\{w\}$.
\item Let $J\subseteq \{1,2,3\}$ of size $2$. If $\{b_j:j\in J\}\subseteq S_4$ and $\{a_j:j\in J\}\subseteq S_i$ for some $i=1,2,3$, then the vertices of $\{b_j:j\in \{1,2,3,4\}\setminus J\}$ are contained in $S_i\cup S_4\cup\{w\}$.
\end{enumerate}
Suppose (C) does not hold. By symmetry, we assume that $J=\{1,2\}$. As $\{a_1, a_2\}$ is not complete to $\{b_1, b_2\}$ in $\bw_3$, we have $i=2$. 
Thus, $a_3\in S_1\cup S_3$. Then $a_1a_2b_2b_1a_1$ is an induced cycle dominated by~$a_3$ in~$G$, and $G$ contains an induced subgraph isomorphic to $W_4$.
A similar argument shows (D) as well.

\begin{enumerate}[label=\rm(\Alph*)]
\setcounter{enumi}{4}
\item If $a_1\in S_4$, $a_2\in S_1$, and $b_1,b_4\in S_2$, then $a_3\notin S_3$.
\end{enumerate}
If $a_3\in S_3$, then $a_1a_2b_4a_3a_1$ is an induced cycle of length $4$ in $G$ dominated by $b_1$. Thus, $G$ contains an induced subgraph isomorphic to $W_4$.

\begin{enumerate}[label=\rm(\Alph*)]
\setcounter{enumi}{5}
\item If $b_1\in S_4$, $b_3\in S_1$, and $a_1,a_2\in S_2$, then $b_4\notin S_3$.
\item If $b_1\in S_4$, $b_3\in S_3$, and $a_1,a_2\in S_2$, then $b_4\notin S_1$.
\end{enumerate}
If $b_4\in S_3$, then $b_1b_3a_2b_4b_1$ is an induced cycle in $G$ dominated by $a_1$. Thus, $G$ contains an induced subgraph isomorphic to $W_4$.
A similar argument shows (G) as well.

\begin{enumerate}[label=\rm(\Alph*)]
\setcounter{enumi}{7}
\item If $a_1,a_2\in S_1$ and $b_1,b_2\in S_3$, then $a_3,b_3,b_4\notin S_1\cup S_3$.
\end{enumerate}
If $a_3\in S_1$, then $a_1b_2b_1a_2a_1$ is an induced cycle in $G$ dominated by $a_3$. Thus, $G$ contains an induced subgraph isomorphic to $W_4$.  If $a_3\in S_3$, then $\{a_3, a_1, b_1\}$ forms $3P_1$ in $G$. Thus, $a_3\notin S_1\cup S_3$. A similar argument shows $b_3, b_4\notin S_1\cup S_3$.

\medskip
\noindent
\textbf{(Case 1. $(G\pivot vw)-v= H$.)}

We claim that~$S_2$ is a clique.
Suppose not.

First assume that $w\in U_1$. Then we may assume that $a_1=w$ and $a_3,b_1\in S_2$ by symmetry.
By Lemma~\ref{lem:easy2}\ref{lem:easy2-ii}, we have $b_3\in S_4$.
Then $a_2\in S_2$ because $S_1\cup S_4$ and $S_3\cup S_4$ are cliques.
By Lemma~\ref{lem:easy2}\ref{lem:easy2-ii}, $b_2\in S_4$.
By (D)  with $J=\{2,3\}$, $b_4\in S_2\cup S_4$, contradicting $S_1\cup S_3\neq \emptyset$.

Now, assume $w\in U_2$. Then we may assume that $b_1=w$.
If $a_1\in S_2$ and $b_3\in S_2$ are non-adjacent, then by Lemma~\ref{lem:easy2}\ref{lem:easy2-ii}, $a_3\in S_4$.
Then $b_2,b_3\in S_2$ because $S_1\cup S_4$ is a clique.
Since $S_1\cup S_3\neq \emptyset$, we have $a_2\in S_3$, contradicting Lemma~\ref{lem:easy2}\ref{lem:easy2-ii} because $a_1, b_2\in S_2$ and $a_2\in S_3$.
This completes the claim.

\medskip
Suppose $w=a_1$.

First assume that $b_3\in S_4$.
Since $S_1\cup S_4$ is a clique, $a_2\in S_2$, and this implies that $b_1\in S_1$ because~$S_2$ is a clique. Again as $S_1$ is a clique, we have $a_3\in S_2$.
By (D) with $J=\{2,3\}$, we have $b_2\in S_3$.
Since~$b_4$ is not a neighbour of~$w$, $b_4$ is in~$S_3\cup S_4$, and 
by (F), $b_4\in S_4$.
Then $a_3b_2b_4b_1a_3$ is an induced cycle of $G$ dominated by $b_3$, and $G$ contains an induced subgraph isomorphic to~$W_4$.
Thus, we may assume by symmetry that $b_2,b_3\notin S_4$. 
This implies that $b_2,b_3\in S_3$.

If $b_4\in S_4$, then $a_2,a_3\in S_2$ because $S_1\cup S_4$ is a clique.
This contradicts Lemma~\ref{lem:easy2}\ref{lem:easy2-v} for the induced cycle $a_2a_3b_3b_2a_2$.
So we may assume that $b_4\in S_3$.
By Lemma~\ref{lem:easy2}\ref{lem:easy2-v}, one of $a_2$ and $a_3$ is not contained in $S_2$.
If both $a_2$ and $a_3$ are contained in $S_1$, then $a_2b_3b_2a_3a_2$ is an induced cycle of $G$ dominated by $b_4$. Thus, we may assume that $a_2$ and $a_3$ are contained in distinct sets of $S_1$ and $S_2$. Since each of $S_1$ and $S_2$ is a clique, $b_1\notin S_1\cup S_2$, a contradiction.

Therefore, we may assume that $w\in U_2$.

\medskip
Suppose $w=b_1$.

First assume that $a_2\in S_4$.
Since $S_1\cup S_4$ is a clique, we have $b_3,b_4\in S_2$.
Since $S_1$ and $S_2$ are cliques, $a_1\in S_1$ and $b_2\in S_2$.
By (C) with $J=\{2,3\}$, $a_3\notin S_4$ and so $a_3\in S_3$.
This contradicts (E).
Thus, we may assume that $a_2,a_3\in S_3$.

By Lemma~\ref{lem:easy2}\ref{lem:easy2-v}, one of $b_2$ and $b_3$ is contained in $S_1$. Assume that $b_2, b_3\in S_1$. Then since $S_1$ and $S_2$ are cliques, we have $a_1\in S_2$ and $b_4\in S_1$. This implies that $b_2a_3a_2b_3b_2$ is an induced cycle of $G$ dominated by $b_4$. Thus, we may assume that $b_2$ and $b_3$ are contained in distinct sets of $S_1$ and $S_2$. Since each of $S_1$ and $S_2$ is a clique, $a_1\notin S_1\cup S_2$, a contradiction.

\medskip
\noindent
\textbf{(Case 2. $(G\pivot vw)-v-w=H$.)}

Suppose $a_1\in S_4$.
Then $b_2,b_3,b_4\in S_2$ because $S_1\cup S_4$ and $S_3\cup S_4$ are cliques.
By (B) with $J=\{2,3\}$, we have $b_1\notin S_1\cup S_3$.
Since $S_1\cup S_3\neq\emptyset$, by symmetry, we may assume that $S_1\neq\emptyset$ and $a_2\in S_1$.
By Lemma~\ref{lem:easy2}\ref{lem:easy2-v} for the induced cycle $a_2b_2b_3a_3a_2$, we have $a_3\in S_3\cup S_4$. Since $a_2\in S_1$ and $S_1$ is complete to $\{w\}\cup S_4$ in $G\pivot vw$, $b_1\notin \{w\}\cup S_4$ and thus $b_1\in S_2$.
By (E), $a_3\notin S_3$ and so $a_3\in S_4$, contradicting (C).
Thus we may assume that $a_1,a_2,a_3\notin S_4$ by symmetry.

\medskip
Suppose $b_1\in S_4$.
Then $a_2,a_3\in S_2$ because $S_1\cup S_4$ and $S_3\cup S_4$ are cliques.
Then by (A) with $J=\{2,3\}$, we have $a_1\notin S_1\cup S_3$ and so $a_1\in S_2\cup S_4$.

If $b_2\in S_1$, then $b_3\in S_3\cup S_4$ by Lemma~\ref{lem:easy2}\ref{lem:easy2-v}. By (D) with $J=\{1,3\}$, we have $b_3\notin S_4$ and so $b_3\in S_3$.
By Lemma~\ref{lem:easy2}\ref{lem:easy2-i} for $\{b_2,b_3,b_4\}$, $b_4\notin S_2$.
By (F) with $J=\{2,3\}$, $b_4\notin S_3$, and 
by (G) with $J=\{2,3\}$, $b_4\notin S_1$ and therefore $b_4\in S_4$.
Then $a_1\in S_2$ and $a_1b_2b_4b_3a_1$ is an induced cycle of~$G$ dominated by $b_1$, a contradiction.

Thus we may assume that $b_2,b_3\notin S_1$ and $b_2,b_3\notin S_3$ by symmetry.
So $b_2,b_3\in S_2\cup S_4$.
Because of Lemma~\ref{lem:easy2}\ref{lem:easy2-v}, one of $b_2$ and $b_3$ is not contained in $S_2$ and so it is contained in $S_4$. This implies that $a_1$ is not contained in $S_1\cup S_3\cup S_4$ and so $a_1\in S_2$.
As $S_1\cup S_3\neq\emptyset$, we have $S_1\cup S_3=\{b_4\}$.
If $b_2\in S_4$, then the induced cycle $a_1a_2b_2b_1a_1$ with $b_4$ contradicts (D). Similarly, if $b_3\in S_4$, then the induced cycle $a_1a_3b_3b_1a_1$ with $b_4$ contradicts (D).

Therefore we may assume that $b_1,b_2,b_3\notin S_4$ by symmetry.

\medskip
Suppose $b_4\in S_4$.
Since $S_1\cup S_4$ and $S_3\cup S_4$ are cliques, it follows that $a_1,a_2,a_3\in S_2$.
Without loss of generality, we may assume that $b_1\in S_1\cup S_3$ because $S_1\cup S_3\neq\emptyset$.
Since $b_2, b_3\notin S_4$, by Lemma~\ref{lem:easy2}\ref{lem:easy2-v}, $b_2,b_3$ are contained in the set of $\{S_1, S_3\}$ different from the set containing $b_1$.
Then $\{a_2,a_3,b_3,b_2\}\subseteq S_2\cup S_3$ induces a cycle in~$G\pivot vw$,
contradicting Lemma~\ref{lem:easy2}\ref{lem:easy2-v}.
Therefore we may assume that $S_4=\emptyset$.

\medskip
We claim that~$S_2$ is a clique.
Suppose not.
As $U_2$ is a clique, $U_1$ must contain a vertex of~$S_2$.
Without loss of generality, we assume that $a_1\in S_2$ and it has a non-neighbour in~$S_2$.
If $b_2\in S_2$, then by Lemma~\ref{lem:easy2}\ref{lem:easy2-ii}, $a_2,b_1\in S_2$.
But this contradicts Lemma~\ref{lem:easy2}\ref{lem:easy2-v}.
So we may assume that $b_2,b_3\notin S_2$ and $b_4\in S_2$.
By Lemma~\ref{lem:easy2}\ref{lem:easy2-ii}, $b_2$ and $b_3$ are contained in the same set of $\{S_1, S_3\}$. Without loss of generality, we assume that $b_2, b_3\in S_1$. By Lemma~\ref{lem:easy2}\ref{lem:easy2-ii} again, we have $b_1\in S_1\cup S_2$.

By Lemma~\ref{lem:easy2}\ref{lem:easy2-v}, $a_2, a_3\in S_3$. If $b_1\in S_1$, then $b_2a_3a_2b_3b_2$ is an induced cycle of $G$ dominated by $b_1$, and $G$ contains $W_4$ as an induced subgraph. If $b_1\in S_2$, then $\{a_1,a_2,a_3,b_1,b_2,b_3,b_4\}$ induces~$\bw_3$ in~$G$.
Therefore we may assume that~$S_2$ is a clique.

\medskip
Then~$V(H)$ is partitioned into three cliques and there exist distinct $i,j\in\{1,2,3\}$ such that $a_i,a_j$ are in one clique and $b_i,b_j$ are in one clique.
We may assume $i=1$ and $j=2$.
By Lemma~\ref{lem:easy2}\ref{lem:easy2-v}, we may assume that $a_1,a_2\in S_1$ and $b_1,b_2\in S_3$.
By (H), $a_3,b_4\in S_2$, contradicting the assumption that~$S_2$ is a clique.
This completes the proof for the case that $(G\pivot vw)-v-w=H$.
\end{proof}

We are now ready to prove the first statement of Theorem~\ref{t-claw}.

\begin{proposition}\label{p-3P_1}
The following statements are equivalent for every graph~$G$:
\begin{enumerate}[label=\rm(\roman*)]
\item $G$ is $3P_1$-pivot-minor-free.
\item $G$ is $(3P_1, W_4, \bw_3)$-free.
(see \figurename~\ref{fig:3P1etc-graphs}).
\end{enumerate}
\end{proposition}
\begin{proof}
We first prove that (i) implies (ii).
Suppose~$G$ contains an induced subgraph $H$ isomorphic to a graph in $\{3P_1,\allowbreak  W_4, \bw_3\}$.
If~$H$ is isomorphic to $W_4$, then by pivoting an edge incident to the vertex of degree~$4$ we obtain a graph which contains an induced subgraph isomorphic to~$3P_1$.
If $H$ is isomorphic to $\bw_3$, then let $U_1=\{a_1,a_2,a_3\}$ and $U_2=\{b_1,b_2,b_3,b_4\}$ be the two cliques of~$H$ and $a_ib_i\in E(H)$ for $i=1,2,3$.
By pivoting an edge~$a_1b_1$, we obtain a subgraph of $G$ induced by $\{a_2,a_3,b_2,b_3,b_4\}$ that is isomorphic to~$W_4$.

Next, we prove that (ii) implies (i).
Suppose~$G$ contains a pivot-minor isomorphic to~$3P_1$.
We use induction on $|V(G)|=n$ to prove that~$G$ contains an induced subgraph isomorphic to a graph in $\{3P_1, W_4, \bw_3\}$.
We may assume that $n \geq 4$.

As $n\geq 4 > |V(3P_1)|$, Lemma~\ref{lem:bouchet} implies that there is a vertex $v\in V(G)$ such that $G-v$ or $(G\pivot vw)-v$, for some neighbour~$w$ of~$v$, contains a pivot-minor isomorphic to~$3P_1$.

If $G-v$ contains a pivot-minor isomorphic to~$3P_1$, then by the induction hypothesis, $G-v$ contains an induced subgraph isomorphic to a graph in $\{3P_1, W_4, \bw_3\}$, hence so does~$G$.
Now we assume that $(G\pivot vw)-v$, for some neighbour~$w$ of~$v$, contains a pivot-minor isomorphic to~$3P_1$.
By the induction hypothesis, $(G\pivot vw)-v$ contains an induced subgraph isomorphic to $3P_1$, $W_4$ or~$\bw_3$.
Applying Lemmas~\ref{lem:subcase1},~\ref{lem:subcase2} and~\ref{lem:subcase3}, respectively, we find that~$G$ contains an induced subgraph isomorphic to a graph in $\{3P_1,W_4,\bw_3\}$.
\end{proof}

As mentioned, for the case where $H=\mbox{claw}$, we need a lemma that allows us to focus on {\em connected} graphs.

\begin{lemma}\label{lem:3P_1-bullclawP5-correspondence}
A graph~$G$ is $(\text{\textnormal{bull}},\text{\textnormal{claw}},P_5)$-free if and only if every connected component of~$G$ is $3P_1$-free.
\end{lemma}

\begin{proof}
The bull, the claw and~$P_5$ are all connected graphs that contain an induced subgraph isomorphic to~$3P_1$.
Therefore, if every connected component of~$G$ is $3P_1$-free, then every connected component of~$G$ is $(\mbox{bull},\mbox{claw},P_5)$-free and so~$G$ is $(\mbox{bull},\mbox{claw},P_5)$-free.

Now suppose that~$G$ contains a connected component containing an induced subgraph isomorphic to~$3P_1$, say on the vertex set $\{x,y,z\}$.
We will show that~$G$ contains an induced subgraph isomorphic to the bull, the claw or~$P_5$.
Let $P=p_1p_2 \cdots p_m$ be a shortest path from $x=p_1$ to $y=p_m$.
If $m\geq 5$, then~$G$ contains an induced subgraph isomorphic to~$P_5$.
Thus, we may assume $m\leq 4$, which in particular implies that~$P$ does not contain~$z$.
Let $Q=q_1q_2 \cdots q_n$ be a shortest path from $z=q_1$ to a vertex~$q_n$ on~$P$.
If~$z$ is adjacent to a vertex on~$P$, then~$G$ must contain the claw or the bull because~$z$ is non-adjacent to~$x$ and~$y$.
So we may assume $n\ge 3$.

If~$q_{n-1}$ has only~$q_n$ as a neighbour on~$P$, then $G[V(P)\cup V(Q)]$ contains an induced subgraph isomorphic to~$P_5$ or the claw.
If~$q_{n-1}$ has exactly two neighbours on~$P$ that are adjacent, then $G[V(P)\cup V(Q)]$ contains an induced subgraph isomorphic to the bull.
In the remaining case, $q_{n-1}$ has two non-adjacent neighbours on~$P$, which means that $q_{n-1}\neq z$ and $G[V(P)\cup V(Q)]$ contains the claw as an induced subgraph.
This completes the proof.
\end{proof}

Combining Lemma~\ref{lem:3P_1-bullclawP5-correspondence} with Proposition~\ref{p-3P_1}, it is easy to prove the second statement of Theorem~\ref{t-claw}.

\begin{proposition}\label{p-claw}
The following statements are equivalent for every graph~$G$:
\begin{enumerate}[label=\rm(\roman*)]
\item $G$ is claw-pivot-minor-free.
\item $G$ is $(\text{\textnormal{claw}}, P_5, \text{\textnormal{bull}}, W_4, \bw_3)$-free.
\end{enumerate}
\end{proposition}
\begin{proof}
We first prove that (i) implies (ii).
Suppose~$G$ contains an induced subgraph $H$ isomorphic to a graph in $\{\mbox{claw}, P_5,\allowbreak \mbox{bull},\allowbreak W_4, \bw_3\}$.
If~$H$ is isomorphic to the claw, then trivially~$G$ contains a pivot-minor isomorphic to the claw.
If~$H$ is isomorphic to $P_5$, then by pivoting the edge between the second and third vertex we obtain a graph which contains an induced claw.
If~$H$ is isomorphic to the bull, then by pivoting an edge incident to the vertex of degree~$2$ we obtain a graph which contains an induced claw.
If~$H$ is isomorphic to $W_4$, then by pivoting an edge incident to the vertex of degree~$4$ we obtain a graph which contains an induced claw.
As shown in the proof of Proposition~\ref{p-3P_1}, pivoting one edge of~$\bw_3$ yields a graph containing an induced~$W_4$.
Thus~$\bw_3$ also contains a pivot-minor isomorphic to the claw.

\smallskip
\noindent
Next, we prove that (ii) implies (i).
Suppose~$G$ contains a pivot-minor isomorphic to the claw.
As edge pivots do not change the number of connected components in a graph and the claw is a connected graph, it follows that~$G$ contains a connected component~$D$ that contains a pivot-minor isomorphic to the claw.

If~$D$ is $(3P_1, W_4,\bw_3)$-free, then by Proposition~\ref{p-3P_1}, $D$ is $3P_1$-pivot-minor-free. This contradicts the fact that $D$ contains a pivot-minor isomorphic to $3P_1$, which is in the claw. Thus, we may assume that $D$ contains an induced subgraph isomorphic to~$3P_1$.
By Lemma~\ref{lem:3P_1-bullclawP5-correspondence}, it follows that~$D$ and therefore~$G$ contains an induced claw, $P_5$ or bull.
This completes the proof.
\end{proof}

\section{Discussion}\label{s-con}

We aim to continue determining the complexity of $H$-{\sc Pivot-Minor}.
We do not know yet if there is a graph~$H$ for which $H$-{\sc Pivot-Minor} is \NP-complete. Interesting open cases are when $H=K_4$ and $H=C_3+P_1$. For both cases we showed that ${\cal F}_H$ contains infinitely many non-isomorphic graphs. Our current techniques for proving polynomial-time solvability of {\sc $H$-Pivot-Minor} is either to prove that the class of $H$-pivot-minor-free graphs has bounded rank-width or else to prove that all graphs in the set~${\cal F}_H$ of minimal forbidden induced subgraphs has bounded size or some polynomial-time verifiable structure. In the latter case we even obtain a polynomial-time algorithm that is certifying.
It would be interesting to know if ${\cal F}_H$ consists of infinitely many non-isomorphic graphs whenever $H$ is a graph and $H'$ is an induced subgraph $H'$, for which ${\cal F}_{H'}$ consists of infinitely many non-isomorphic graphs.

Besides the above,
a proof for the Minor Recognition conjecture~\cite{GGW06} for binary matroids would also yield a technique to obtain  complexity results for pivot-minors.
In particular, if this conjecture is true, then for every graph~$H$ the {\sc $H$-Pivot-Minor} problem is polynomial-time solvable for bipartite graphs.
This follows from Lemma~\ref{lem:bin}, which implies that a bipartite connected graph~$H$ is a pivot-minor of a bipartite graph~$G$ if and only if for binary matroids~$M$ and~$N$ that have~$G$ and~$H$ as fundamental graphs, respectively, $N$ or the dual of~$N$ is a minor of~$M$ (if~$H$ is not connected, then we try all possible ways of making duals per connected component of~$H$).

\subsection*{Declarations}\ \\
\noindent
\textbf{Data Availability Statement.}
Data sharing not applicable to this article as no datasets were
generated or analysed during the current study.

\medskip 
\noindent
\textbf{Competing Interests Policy.}
The authors have no competing interests as defined by Springer, or other interests that might be perceived to influence the results and/or discussion reported in this paper.

\medskip 
\noindent
\textbf{Author Contributions.}
All authors contributed to the paper.


\begin{thebibliography}{10}

\bibitem{Bo94}
Andr\'e Bouchet.
\newblock Circle graph obstructions.
\newblock {\em Journal of Combinatorial Theory, Series B}, 60:107--144, 1994.

\bibitem{BD91}
Andr{\'e} Bouchet and Alain Duchamp.
\newblock Representability of {$\Delta$}-matroids over {${\rm GF}(2)$}.
\newblock {\em Linear Algebra and Its Applications}, 146:67--78, 1991.

\bibitem{BV87}
A.~E. Brouwer and H.~J. Veldman.
\newblock Contractibility and {NP}-completeness.
\newblock {\em Journal of Graph Theory}, 11:71--79, 1987.

\bibitem{CSSS20}
Maria Chudnovsky, Alex Scott, Paul~D. Seymour, and Sophie Spirkl.
\newblock Detecting an odd hole.
\newblock {\em Journal of the {A}{C}{M}}, 67:5:1--5:12, 2020.

\bibitem{CLB81}
D.~G. Corneil, H.~Lerchs, and L.~Stewart Burlingham.
\newblock Complement reducible graphs.
\newblock {\em Discrete Applied Mathematics}, 3:163--174, 1981.

\bibitem{CO07}
Bruno Courcelle and Sang-il Oum.
\newblock Vertex-minors, monadic second-order logic, and a conjecture by
  {S}eese.
\newblock {\em Journal of Combinatorial Theory, Series B}, 97:91--126, 2007.

\bibitem{DDJKKOP21}
Konrad~K. Dabrowski, Fran\c{c}ois Dross, Jisu Jeong, Mamadou~M. Kant\'{e},
  O{-}joung Kwon, Sang{-}il Oum, and Daniel Paulusma.
\newblock Tree pivot-minors and linear rank-width.
\newblock {\em SIAM Journal on Discrete Mathematics}, 35:2922--2945, 2021.

\bibitem{Dabrowski2018}
Konrad~K. Dabrowski, Fran\c{c}ois Dross, Jisu Jeong, Mamadou~Moustapha
  Kant\'{e}, O{-}joung Kwon, Sang{-}il Oum, and Dani\"{e}l Paulusma.
\newblock Computing small pivot-minors.
\newblock In {\em Proceedings of the 44th {I}nternational {W}orkshop on
  {G}raph-{T}heoretic {C}oncepts in {C}omputer {S}cience, {W}{G} 2018}, volume
  11159 of {\em Lecture Notes in Computer Science}, pages 125--138. Springer,
  Heidelberg, 2018.

\bibitem{Dahlberg2022}
Axel Dahlberg, Jonas Helsen, and Stephanie Wehner.
\newblock The complexity of the vertex-minor problem.
\newblock {\em Information Processing Letters}, 175:106222, 2022.

\bibitem{FKMP95}
M.~R. Fellows, J.~Kratochv{\'\i}l, M.~Middendorf, and F.~Pfeiffer.
\newblock The complexity of induced minors and related problems.
\newblock {\em Algorithmica}, 13:266--282, 1995.

\bibitem{GJT76}
M.~R. Garey, D.~S. Johnson, and R.~Endre Tarjan.
\newblock The planar {H}amiltonian circuit problem is {NP}-complete.
\newblock {\em SIAM Journal on Computing}, 5:704--714, 1976.

\bibitem{GGW06}
Jim Geelen, Bert Gerards, and Geoff Whittle.
\newblock Towards a structure theory for matrices and matroids.
\newblock In {\em International {C}ongress of {M}athematicians. {V}ol. {III}},
  pages 827--842. European Mathematical Society, Z\"urich, 2006.

\bibitem{GKMW19}
Jim Geelen, O{-}joung Kwon, Rose McCarty, and Paul Wollan.
\newblock The grid theorem for vertex-minors.
\newblock {\em Journal of Combinatorial Theory, Series B}, 158:93--16, 2023.

\bibitem{GO09}
Jim Geelen and Sang-il Oum.
\newblock Circle graph obstructions under pivoting.
\newblock {\em Journal of Graph Theory}, 61:1--11, 2009.

\bibitem{GKMW11}
Martin Grohe, Ken-ichi Kawarabayashi, D\'aniel Marx, and Paul Wollan.
\newblock Finding topological subgraphs is fixed-parameter tractable.
\newblock In {\em {P}roceedings of the 43rd {ACM} {S}ymposium on {T}heory of
  {C}omputing, {S}{T}{O}{C} 2011}, pages 479--488. ACM, New York, 2011.

\bibitem{KK18}
Mamadou~Moustapha Kant\'{e} and O{-}joung Kwon.
\newblock Linear rank-width of distance-hereditary graphs {II}. {V}ertex-minor
  obstructions.
\newblock {\em European Journal of Combinatorics}, 74:110--139, 2018.

\bibitem{KMOW21}
O{-}joung Kwon, Rose McCarty, Sang{-}il Oum, and Paul Wollan.
\newblock Obstructions for bounded shrub-depth and rank-depth.
\newblock {\em Journal of Combinatorial Theory, Series B}, 149:76--91, 2021.

\bibitem{KO14}
O-joung Kwon and Sang-il Oum.
\newblock Graphs of small rank-width are pivot-minors of graphs of small
  tree-width.
\newblock {\em Discrete Applied Mathematics}, 168:108--118, 2014.

\bibitem{LLMT09}
Benjamin L\'ev\^eque, David~Y. Lin, Fr\'ed\'eric Maffray, and Nicolas
  Trotignon.
\newblock Detecting induced subgraphs.
\newblock {\em Discrete Applied Mathematics}, 157:3540--3551, 2009.

\bibitem{LMT12}
Benjamin L\'ev\^eque, Fr\'ed\'eric Maffray, and Nicolas Trotignon.
\newblock On graphs with no induced subdivision of {$K_4$}.
\newblock {\em Journal of Combinatorial Theory, Series B}, 102:924--947, 2012.

\bibitem{MT92}
Ji{\v r}{\'\i} Matou{\v s}ek and Robin Thomas.
\newblock On the complexity of finding iso- and other morphisms for partial
  {$k$}-trees.
\newblock {\em Discrete Mathematics}, 108:343--364, 1992.

\bibitem{MMNS11}
Ross~M. Mcconnell, Kurt Mehlhorn, Stefan N\"{a}her, and Pascal Schweitzer.
\newblock Survey: Certifying algorithms.
\newblock {\em Computer Science Review}, 5:119--161, 2011.

\bibitem{Ou05}
Sang-il Oum.
\newblock Rank-width and vertex-minors.
\newblock {\em Journal of Combinatorial Theory, Series B}, 95:79--100, 2005.

\bibitem{Ou08}
Sang-il Oum.
\newblock Rank-width and well-quasi-ordering.
\newblock {\em SIAM Journal on Discrete Mathematics}, 22:666--682, 2008.

\bibitem{Ou09}
Sang-il Oum.
\newblock Excluding a bipartite circle graph from line graphs.
\newblock {\em Journal of Graph Theory}, 60:183--203, 2009.

\bibitem{Ou16}
Sang-il Oum.
\newblock Rank-width: Algorithmic and structural results.
\newblock {\em Discrete Applied Mathematics}, 231:15--24, 2017.

\bibitem{OS06}
Sang-il Oum and Paul Seymour.
\newblock Approximating clique-width and branch-width.
\newblock {\em Journal of Combinatorial Theory, Series B}, 96:514--528, 2006.

\bibitem{Oxley}
James Oxley.
\newblock {\em Matroid {T}heory}, volume~21 of {\em Oxford Graduate Texts in
  Mathematics}.
\newblock Oxford University Press, Oxford, second edition, 2011.

\bibitem{RS95}
Neil Robertson and P.~D. Seymour.
\newblock Graph minors. {XIII}. {T}he disjoint paths problem.
\newblock {\em Journal of Combinatorial Theory, Series B}, 63:65--110, 1995.

\bibitem{HKPST12}
Pim van~'t Hof, Marcin Kami\'nski, Dani\"el Paulusma, Stefan Szeider, and
  Dimitrios~M. Thilikos.
\newblock On graph contractions and induced minors.
\newblock {\em Discrete Applied Mathematics}, 160:799--809, 2012.

\end{thebibliography}
\end{document}